\documentclass[12pt]{article}
\usepackage{amssymb}
\newcommand{\R}{\mathbb{R}}

\newcommand{\Z}{\mathbb{Z}}
\renewcommand{\P}{\mathbb{P}}
\newcommand{\Q}{\mathbb{Q}}
\newcommand{\E}{\mathbb{E}}

\newcommand{\N}{\mathbb{N}}

\newcommand{\T}{\mathbb{T}}
\def\build#1_#2^#3{\mathrel{
\mathop{\kern 0pt#1}\limits_{#2}^{#3}}}
\def\llbracket{[\hspace{-.10em} [ }
\def\rrbracket{ ] \hspace{-.10em}]}
\def\cq{$\hfill \square$}
\def \un{\underline}

\def\g{{\cal G}}

\def\e{{\cal E}}

\def\u{{\cal U}}

\def\t{{\cal T}}
\def\v{{\cal V}}

\def\q{{\cal Q}}

\def\z{{\cal Z}}

\def\eg{{\bf e}}

\def\be{\begin{equation}}
\def\ee{\end{equation}}
\def\ba{\begin{eqnarray*}}
\def\ea{\end{eqnarray*}}
\def\ov{\overline}
\def\wh{\widehat}
\def\wt{\widetilde}

\def\la{\longrightarrow}
\def\da{\downarrow}

\def\noi{\noindent}
\def\proof{\vskip 3mm \noindent{\bf Proof:}\hskip10pt}

\newtheorem{theorem}{Theorem}[section]
\newtheorem{lemma}[theorem]{Lemma}
\newtheorem{proposition}[theorem]{Proposition}
\newtheorem{corollary}[theorem]{Corollary}

\begin{document}

\title{ \bf AN INVARIANCE PRINCIPLE FOR 
CONDITIONED TREES}
\author{
Jean-Fran\c cois {\sc Le Gall}\\
{\small D.M.A., Ecole normale sup\'erieure, 
45 rue d'Ulm, 75005 Paris, France}\\
{\small legall@dma.ens.fr}}
\date{} 
\maketitle


\begin{abstract}
We consider Galton-Watson trees associated with a critical 
offspring distribution and conditioned to have exactly $n$ vertices.
These trees are embedded in the real line by affecting spatial 
positions to the vertices, in such a way that the increments of the
spatial positions along edges of the tree are independent variables
distributed according to a symmetric probability distribution on
the real line. We then condition on the event that all
spatial positions are nonnegative.
Under suitable assumptions on the offspring 
distribution and the spatial displacements, we prove that
these conditioned spatial trees converge as $n\to\infty$, modulo an
appropriate rescaling, towards the conditioned Brownian tree
that was studied in previous work. Applications are given to
asymptotics for random quadrangulations.
\end{abstract}

\section{Introduction}

The main goal of the present work is to prove an
invariance principle for Galton-Watson trees embedded
in the real line and constrained to remain
on the positive side. One major motivation for
this problem came from recent asymptotic results 
for random quadrangulations which have been established
by Chassaing and Schaeffer \cite{CS}.

The asymptotic behavior of Galton-Watson trees conditioned
to have a large fixed progeny was investigated by Aldous \cite{Al3}
in connection with the so-called Continuum Random Tree (CRT). Precisely,
under the assumption that the offspring distribution $\mu$ is critical
and has finite variance $\sigma^2>0$, a Galton-Watson tree conditioned to have 
exactly $n$ vertices, with edges rescaled by the factor $\sigma n^{-1/2}/2$,
will converge in distribution, in a suitable sense, towards the
CRT. A convenient way of making this convergence mathematically precise
is to use the contour function of the conditioned 
Galton-Watson tree (cf Fig.1 below). Modulo a rescaling analogous
to the classical Donsker theorem for random walks, this contour function
converges in distribution as $n\to\infty$ towards a normalized Brownian excursion, that is
a positive Brownian excursion conditioned to have duration $1$ 
(cf the convergence of the first components in Theorem \ref{invar0} below).
Informally we may say that the normalized Brownian excursion is the contour 
function of the CRT. See \cite{DLG} for analogous contour descriptions
of the more general L\'evy trees, and \cite{Duq} for a recent generalization of
Aldous' theorem.

In view of various applications, and in particular in connection
with the theory of superprocesses, it is interesting to combine 
the branching structure of the Galton-Watson tree with spatial displacements.
Here we consider the simple special case where these
spatial displacements are given by a one-dimensional symmetric random walk
on the real line with jump distribution
$\gamma$. This means that i.i.d.
random variables $Y_e$ with distribution
$\gamma$ are associated with the different edges of the tree, and that the spatial
position $U_v$ of a vertex $v$ is obtained by summing the displacements $Y_e$ corresponding
to edges $e$ that belong to the path from the root to the vertex $v$. The resulting object,
called here a spatial tree, consists of a (random) pair $(\t,U)$, where $\t$ is a 
discrete (plane) tree and $U$ is a mapping from the set of vertices of $\t$ into $\R$.
In the same way as the tree $\t$ can be coded by its contour function, a convenient
way of encoding the spatial positions is via the spatial contour function
(see Section 2 for a precise definition, and Fig.2 for an example).

From now on, we suppose that the tree $\t$ is a Galton-Watson tree with offspring 
distribution $\mu$ satisfying the above assumptions, and conditioned to have 
exactly $n$ vertices, and that the spatial positions $U_v$ are generated as explained
in the preceding paragraph.
We assume furthermore that $\mu$ has (small) exponential moments and
that $\gamma([x,\infty))=o(x^{-4})$ as $x\to\infty$. We denote by
$\rho^2$ the variance of $\gamma$. Then rescaling both the edges of $\t$ by the
factor $\sigma n^{-1/2}/2$ and the spatial displacements $U_v$ by $\rho^{-1}(\sigma/2)^{1/2}n^{-1/4}$
will lead as $n\to\infty$ to a limiting object which is independent of $\mu$ and $\gamma$.
A precise statement for this convergence is given in Theorem \ref{invar0} below, which is taken
from Janson and Marckert \cite{JM} (see \cite{CS}, \cite{Git} and \cite{MM2}
for similar statements under different assumptions -- related results have also been
obtained by Kesten \cite{Ke} under other conditionings of the tree). 
This convergence is closely related to the approximation of superprocesses
by branching particle systems: See in particular \cite{LG0}. Roughly
speaking, the limiting object combines the branching structure of the CRT with spatial displacements
given by independent linear Brownian motions along the edges of the tree. A
convenient representation of this limiting object,
which is used in Theorem \ref{invar0}, is provided by the Brownian snake (see e.g. \cite{Zu}).
To describe this approach, let $r\in\R$,
which will represent the initial position (the spatial position of the root) and let
$\eg=(\eg(s),0\leq s\leq 1)$ be a normalized Brownian excursion. Let $Z^r=(Z^r(s),0\leq s\leq 1)$ be
a real-valued process such that, conditionally given
$\eg$, $Z^r$ is Gaussian with mean and covariance given by
\begin{description}
\item{$\bullet$} $E[Z^r(s)]=r$ for every $s\in[0,1]$;
\item{$\bullet$} ${\rm cov}(Z^r(s),Z^r(s'))=\displaystyle{\inf_{s\leq t\leq s'} \eg(t)}$
for every $0\leq s\leq s'\leq 1$.
\end{description}
Informally, each time $s\in[0,1]$ corresponds via the contour function coding to a vertex of the CRT,
and $Z^r(s)$ is the spatial position of this vertex. The formula for the conditional
covariance of $Z^r(s)$ and $Z^r(s')$ is then justified by the fact that $\inf_{s\leq t\leq s'}
\eg(t)$ is the generation of the ``most recent common ancestor'' to the vertices 
corresponding to $s$ and $s'$  (see the
introduction to
\cite{LGW} for a precise version of this informal explanation). In the terminology of
\cite{Zu}, $Z^r$ is the terminal point process of the one-dimensional Brownian snake driven by
the normalized Brownian excursion $\eg$ and with initial point $r$.
For simplicity, we will say here that $(\eg,Z^r)$ is the Brownian snake with initial point $r$. The
random measure known as one-dimensional ISE (Integrated Super-Brownian Excursion) may be defined by
the formula:
$$\langle\z,f\rangle=\int_0^1 ds\,f(Z^0_s)\;,\quad f\in C_b(\R).$$
ISE in higher dimensions has found many applications in asymptotics for models of 
statistical physics: See in particular \cite{Sl1}, \cite{Sl2} and \cite{Sl3}.

Our main interest in this work is to study asymptotics for the discrete spatial trees $(\t,U)$
conditioned on the event that spatial positions $U_v$ remain in the positive half-line. In contrast
with the situation described above, this induces an interaction between the branching structure and
the spatial displacements, which makes the analysis of the model more delicate. We start from a
spatial tree
$(\t,U)$ generated as explained above. Then our main result
(Theorem \ref{invar}) states that this spatial tree 
conditioned on the event that $U_v\geq 0$ for every vertex $v$ of $\t$
and rescaled as previously
will converge in distribution as $n\to\infty$ to a (universal) limiting object.
This limiting object corresponds to the conditioned Brownian snake that was
studied in detail in \cite{LGW}. More precisely, for every $r>0$, let $(\ov \eg^r,\ov Z^r)$
be distributed as the pair $(\eg,Z^r)$ introduced above, under the conditioning that
$Z^r(s)\geq 0$ for every $s\in[0,1]$. Note that as long as $r>0$, this conditioning
involves an event of positive probability. Theorem 1.1 in \cite{LGW} shows that 
the process $(\ov \eg^r,\ov Z^r)$ converges in distribution as $r\da 0$ towards a
limiting pair $(\ov \eg^0,\ov Z^0)$, which is our conditioned Brownian snake
with initial point $0$.
According to Theorem \ref{invar} below, the pair $(\ov \eg^0,\ov Z^0)$ is the scaling limit of 
the pair consisting of the contour function and the spatial contour function of our
spatial trees conditioned to have nonnegative spatial positions.

The preceding description of the conditioned process $(\ov \eg^0,\ov Z^0)$ 
as the limit of $(\ov \eg^r,\ov Z^r)$ when $r\downarrow 0$
does not give much information about this process. Note in particular that the underlying conditioning
is in a sense very degenerate, since we are dealing with a continuous tree 
of linear Brownian paths all started from the origin and conditioned not to hit
the negative half-line. Still Theorem 1.2 in \cite{LGW} provides a useful construction
of the conditioned object $(\ov \eg^0,\ov Z^0)$ from the unconditioned one 
$(\eg,Z^0)$. In order to present this construction, first recall that there is
a.s. a unique $s_*$ in $(0,1)$ such that
$$Z^0(s_*)=\inf_{s\in[0,1]} Z^0(s)$$
(see Lemma 16 in \cite{MM} or Proposition 2.5 in \cite{LGW}). For every $s\in[0,\infty)$,
write $\{s\}$ for the fractional part of $s$. 
According to Theorem 1.2 in \cite{LGW}, the conditioned snake $(\ov \eg^0,\ov Z^0)$
may be be constructed explicitly as follows: For every $s\in[0,1]$,
\ba
&&\ov \eg^0(s) = \eg({s_*}) + \eg(\{s_*+s\}) - 
2 \,\inf_{s\wedge\{s_*+s\}\leq u\leq s\vee \{s_*+s\}} \eg(u)\\
&&\ov Z^0(s) = Z^0(\{s_*+s\}) - Z^0({s_*})\,.
\ea
In terms of trees, this means that the conditioned tree is obtained by re-rooting
the unconditioned one at the vertex having the minimal spatial position: 
The above formula for $\ov \eg^0(s)$ corresponds to the contour function for the
tree coded by $\eg$ re-rooted at the vertex $s_*$ -- see the discussion in the
introduction of \cite{LGW}. This construction of the conditioned snake $(\ov \eg^0,\ov Z^0)$
is of course reminiscent of a famous result of Verwaat 
\cite{Verwaat} connecting the normalized
Brownian excursion and the Brownian bridge.

The initial motivation for the present work came from a recent paper of
Chassaing and Schaeffer \cite{CS}
discussing asymptotics for planar maps (see also \cite{CD} for 
a related work). A key result (Theorem 1 in
\cite{CS} or Theorem \ref{bijection} below) establishes a
bijection between rooted quadrangulations 
with $n$ faces and the so-called well-labelled
trees with $n$ edges.  In the terminology of the
present work, a well-labelled tree is a spatial tree $(\t,U)$ with the additional
properties that the spatial positions are positive integers, the spatial position of the root
is $1$ and the spatial positions of two neighboring vertices can differ by at most $1$. 
The preceding bijection between quadrangulations and trees has the nice
feature that distances of vertices of the quadrangulation from the root
correspond to spatial positions in the associated tree. This suggests that asymptotic
properties of distances from the root in
random quadrangulations with $n$ faces can be read from asymptotics for well-labelled
trees with $n$ edges, an idea which was exploited in \cite{CS}. As an application
of Theorem \ref{invar}, we provide a direct proof of some of the main results of \cite{CS}.
The key idea is to observe that uniform well-labelled trees with $n$ edges can be obtained
as conditioned spatial trees generated by letting the offspring distribution be geometric
with parameter $1/2$, and the spatial distribution be uniform on $\{-1,0,1\}$. 
It then follows from Theorem \ref{invar} that the scaling limit of well-labelled trees
with $n$ vertices is described by our conditioned Brownian snake. As a consequence,
several quantities such as the (suitably rescaled) radius of the quadrangulation have a limit in
distribution which can be expressed, first in terms of the conditioned Brownian snake
and then via the Verwaat-like transformation in terms of the unconditioned snake.

Another recent paper \cite{MM} of Marckert and Mokkadem proposes a model called the
Brownian map for the continuous limit of rooted quadrangulations with $n$ faces.
This construction makes a heavy use of the conditioned Brownian snake, which
is defined in \cite{MM} via the Verwaat transformation rather than as a conditioned object as 
here or in \cite{LGW}. Theorem \ref{invar} readily gives a positive answer to a conjecture of
\cite{MM} (cf Remark 6 in
\cite{MM}) asserting that the Brownian map is, in some sense, the scaling limit of uniform rooted
quadrangulations with $n$ faces.

Connections between trees and planar maps are also of interest in theoretical physics: See 
in particular \cite{BDG0}, \cite{BDG1}, \cite{BDG2} and \cite{Dur}
for discussions and various applications. In this perspective, 
quadrangulations, or more general planar maps, serve as a model of random geometry.
We mention the recent article \cite{BDG3}, which extends the bijection 
between quadrangulations and well-labelled trees to more general classes of
planar maps. It seems plausible that the bijections of \cite{BDG3} can be used in
connection with the methods of the present work in order to generalize the
Chassaing-Schaeffer asymptotics to general planar maps. We hope to address this problem
in some future work.

The present paper is organized as follows. Section 2 gives the basic assumptions and states
our main result Theorem \ref{invar}. Section 3 introduces the key technical idea of
comparing the distribution of spatial trees re-rooted at the minimal spatial position with
that of conditioned trees. In the continuous framework, the Verwaat transformation 
shows that the distribution of the re-rooted tree and that of the conditioned one
are identical. This is no longer true in a discrete setting, but there are still
close relations between the two distributions, that play a major role in our proofs.
As a first application, Section 4 derives estimates for the probability 
that the spatial positions are all positive: This probability is bounded above and
below by constants times $n^{-1}$, where $n$ is the number of vertices
in the tree (Proposition \ref{positive}). Section 5
briefly discusses a spatial Markov property for our spatial trees, which holds for the subtrees
originating from the ``first'' vertices whose spatial position exceeds a level $a>0$. Section 6 again
applies the results of Section 3 to derive some asymptotic regularity properties of conditioned trees.
These regularity properties are first established for unconditioned trees via
Theorem \ref{invar0}. They can then be transferred to conditioned
trees thanks to Section 3. Section 7 gives the proof of Theorem \ref{invar}. Roughly speaking,
the argument goes as follows. Thanks to the regularity properties of Section 6
and the spatial Markov property of Section 5, the conditioned tree is well approximated
by a spatial tree (with a number of vertices of order $n$) whose root is located at a point close to
$\alpha n^{1/4}$ (where $\alpha$ is a ``small'' positive number) and which is conditioned not to hit
the negative half-line. The limit theorem for unconditioned spatial trees (Theorem \ref{invar0})
then shows that the limit of this suitably rescaled spatial tree is described by the
Brownian snake $(\ov \eg^\alpha,\ov Z^\alpha)$ with initial point $\alpha$ and conditioned not
to hit the negative half-line. Notice that the conditioning is not degenerate here since
$\alpha>0$. To complete the proof, we just have to use the fact that $(\ov \eg^\alpha,\ov Z^\alpha)$
is close in distribution to $(\ov \eg^0,\ov Z^0)$ when $\alpha$ is small, as was mentioned above.
Finally, Section 8 discusses applications to random quadrangulations.

\section{Basic assumptions and statement of the main result}

Throughout this work, we denote by $\mu$ the underlying 
offspring distribution, which is a probability measure on
$\Z_+=\{0,1,2,\ldots\}$. We always assume that
$\mu(1)<1$ and
\begin{description}
\item{$\bullet$} $\mu$ is critical, that is 
$\displaystyle \sum_{k=0}^\infty k\mu(k)=1$.
\item{$\bullet$} $\mu$ is aperiodic, that is 
$\mu$ is not supported on a proper subgroup of 
$\Z$.
\item{$\bullet$} $\mu$ has exponential moments: There exists 
a constant $\lambda>0$ such that
$$\sum_{k=0}^\infty \mu(k)\,e^{\lambda k} <\infty.$$
\end{description}
We denote by $\sigma^2>0$ the variance of $\mu$.

The law of the spatial displacement is denoted by $\gamma$. Thus 
$\gamma$ is a probability distribution on $\R$. We exclude the trivial
case $\gamma=\delta_0$ and we always assume that $\gamma$ is symmetric,
that is $\gamma$ is invariant under the transformation $x\la -x$. 
Our main results also require the additional assumption
\be
\label{4mom}
\lim_{x\to \infty} x^{4}\,\gamma([x,\infty))=0.
\ee
When (\ref{4mom}) holds, we denote by $\rho^2>0$ the variance of $\gamma$.

Let us now introduce some formalism for discrete trees, which we borrow from Neveu \cite{Neveu}. Let
$${\cal U}=\bigcup_{n=0}^\infty \N^n  $$
where $\N=\{1,2,\ldots\}$ and by convention $\N^0=\{\varnothing\}$. If
$u=(u_1,\ldots u_m)$ and 
$v=(v_1,\ldots, v_n)$ belong to $\cal U$, we write $uv=(u_1,\ldots u_m,v_1,\ldots ,v_n)$
for the concatenation of $u$ and $v$. In particular $u\varnothing=\varnothing u=u$.
If $p\geq 1$ is an integer, we will write $1^p$ for the $p$-tuple $(1,1,\ldots,1)\in\N^p$.
Finally, we set $\u^*=\u\backslash\{\varnothing\}$.

A plane tree $\t$ is a finite subset of
$\cal U$ such that:
\begin{description}
\item{(i)} $\varnothing\in \t$.

\item{(ii)} If $v\in \t$ and $v=uj$ for some $u\in {\cal U}$ and
$j\in\N$, then $u\in\t$.

\item{(iii)} For every $u\in\t$, there exists a number $N_u(\t)\geq 0$
such that, for every $j\in\N$, $uj\in\t$ if and only if $1\leq j\leq N_u(\t)$.
\end{description}

\noi We denote by ${\bf T}$ the set of all plane trees.

\begin{center}
\unitlength=1.1pt
\begin{picture}(200,130)

\thicklines \put(50,0){\line(-1,2){20}}
\thicklines \put(50,0){\line(1,2){20}}
\thicklines \put(30,40){\line(0,1){40}}
\thicklines \put(30,40){\line(-1,2){20}}
\thicklines \put(30,40){\line(1,2){20}}
\thicklines \put(30,80){\line(-1,2){20}}
\thicklines \put(30,80){\line(1,2){20}}

\thinlines \put(37,15){\vector(-1,2){6}}
\thinlines \put(17,55){\vector(-1,2){6}}
\thinlines \put(20,67){\vector(1,-2){6}}
\thinlines \put(27,57){\vector(0,1){12}}
\thinlines \put(17,95){\vector(-1,2){6}}
\thinlines \put(22,107){\vector(1,-2){6}}
\thinlines \put(33,95){\vector(1,2){6}}
\thinlines \put(48,107){\vector(-1,-2){6}}
\thinlines \put(33,69){\vector(0,-1){12}}
\thinlines \put(34,55){\vector(1,2){6}}
\thinlines \put(48,67){\vector(-1,-2){6}}
\thinlines \put(40,27){\vector(1,-2){6}}
\thinlines \put(54,15){\vector(1,2){6}}
\thinlines \put(68,27){\vector(-1,-2){6}}

\put(40,-5){$\varnothing$}
\put(20,35){$1$}
\put(60,35){$2$}
\put(-3,75){$11$}
\put(17,75){$12$}
\put(37,75){$13$}
\put(-7,115){$121$}
\put(31,115){$122$}
\put(30,-25){tree $\t$}
\thinlines \put(100,0){\line(1,0){120}}
\thinlines \put(100,0){\line(0,1){130}}
\thicklines \put(100,0){\line(1,5){8}}
\thicklines \put(108,40){\line(1,5){8}}
\thicklines \put(116,80){\line(1,-5){8}}
\thicklines \put(124,40){\line(1,5){8}}
\thicklines \put(132,80){\line(1,5){8}}
\thicklines \put(140,120){\line(1,-5){8}}
\thicklines \put(148,80){\line(1,5){8}}
\thicklines \put(156,120){\line(1,-5){8}}
\thicklines \put(164,80){\line(1,-5){8}}
\thicklines \put(172,40){\line(1,5){8}}
\thicklines \put(180,80){\line(1,-5){8}}
\thicklines \put(188,40){\line(1,-5){8}}
\thicklines \put(196,0){\line(1,5){8}}
\thicklines \put(204,40){\line(1,-5){8}}

\thinlines \put(108,0){\line(0,1){2}}
\thinlines \put(116,0){\line(0,1){2}}
\thinlines \put(124,0){\line(0,1){2}}
\thinlines \put(100,40){\line(1,0){2}}
\thinlines \put(100,80){\line(1,0){2}}
\thinlines \put(212,0){\line(0,1){2}}

\put(115,-25){contour function $C(t)$}

\put(107,-6){\scriptsize 1}
\put(115,-6){\scriptsize 2}
\put(123,-6){\scriptsize 3}

\put(95,39){\scriptsize 1}
\put(95,79){\scriptsize 2}

\end{picture}

\vskip 10mm

Figure 1
\end{center}

To define now the {\it contour function} of $\t$, consider a particle that starts from the root
and visits continuously  all edges at speed one, going backwards as less as possible and respecting 
the lexicographical order of vertices. Since each edge will
be crossed twice, the total time needed to explore the tree is
$2(|\t|-1)$, where $|\t|$ denotes the cardinality (number of vertices)
of $\t$. For every $t\in[0,2(|\t|-1)]$,
we let
$C(t)$ denote the distance from the root of the position of the particle at time $t$.
Fig.1 explains
the definition of the contour function better than a formal definition. Clearly a tree $\t$
is uniquely determined by its contour function. 

A (discrete) spatial tree is a pair $(\t,U)$, where $\t\in{\bf T}$
and $U=(U_v,v\in\t)$ is a mapping from the set $\t$ into $\R$. We denote by
$\Omega$ the set of all spatial trees. A spatial tree $(\t,U)$ can be coded
by a pair $(C,V)$, where $C=(C(t),0\leq t\leq 2(|\t|-1))$ is the contour function of 
$\t$ and the {\it spatial contour function} $V=(V(t),0\leq t\leq 2(|\t|-1))$ is
defined as follows. First if $t$ is an integer, then $t$ corresponds 
in the evolution of the contour to a vertex $v$ of $\t$, and 
we put $V(t)=U_v$. We then complete the definition of 
$V$ by interpolating linearly between successive integers.
See Fig.2 for an example: The tree on the left side of this 
figure is the same as in Fig.1, but the numbers in bold attached 
to the different vertices now represent the spatial positions.

\begin{center}
\unitlength=1.1pt
\begin{picture}(200,132)

\thicklines \put(50,0){\line(-1,2){20}}
\thicklines \put(50,0){\line(1,2){20}}
\thicklines \put(30,40){\line(0,1){40}}
\thicklines \put(30,40){\line(-1,2){20}}
\thicklines \put(30,40){\line(1,2){20}}
\thicklines \put(30,80){\line(-1,2){20}}
\thicklines \put(30,80){\line(1,2){20}}
\put(40,-5){$\bf 1$}
\put(20,35){$\bf 3$}
\put(52,35){$\bf -1$}
\put(3,75){$\bf 2$}
\put(23,75){$\bf 1$}
\put(41,75){$\bf 3$}
\put(-7,115){$\bf -1$}
\put(39,115){$\bf 0$}

\thinlines\put(108,30){\line(0,1){2}}
\thinlines\put(116,30){\line(0,1){2}}
\thinlines\put(124,30){\line(0,1){2}}
\thinlines\put(132,30){\line(0,1){2}}
\thinlines\put(140,30){\line(0,1){2}}
\thinlines\put(148,30){\line(0,1){2}}
\thinlines\put(156,30){\line(0,1){2}}
\thinlines\put(164,30){\line(0,1){2}}
\thinlines\put(172,30){\line(0,1){2}}
\thinlines\put(180,30){\line(0,1){2}}
\thinlines\put(188,30){\line(0,1){2}}
\thinlines\put(196,30){\line(0,1){2}}
\thinlines\put(204,30){\line(0,1){2}}
\thinlines\put(212,30){\line(0,1){2}}

\thinlines\put(100,54){\line(1,0){2}}
\thinlines\put(100,78){\line(1,0){2}}
\thinlines\put(100,102){\line(1,0){2}}
\thinlines\put(100,6){\line(1,0){2}}

\put(92,52){$1$}
\put(92,76){$2$}
\put(92,100){$3$}
\put(84,4){$-1$}

\thinlines \put(100,30){\line(1,0){120}}
\thinlines \put(100,0){\line(0,1){130}}
\thicklines \put(100,54){\line(1,6){8}}
\thicklines \put(108,102){\line(1,-3){8}}
\thicklines \put(116,78){\line(1,3){8}}
\thicklines \put(124,102){\line(1,-6){8}}
\thicklines \put(132,54){\line(1,-6){8}}
\thicklines \put(140,6){\line(1,6){8}}
\thicklines \put(148,54){\line(1,-3){8}}
\thicklines \put(156,30){\line(1,3){8}}
\thicklines \put(164,54){\line(1,6){8}}
\thicklines \put(172,102){\line(1,0){8}}
\thicklines \put(180,102){\line(1,0){8}}
\thicklines \put(188,102){\line(1,-6){8}}
\thicklines \put(196,54){\line(1,-6){8}}
\thicklines \put(204,6){\line(1,6){8}}

\put(-6,-25){spatial tree $(\t,U)$}

\put(102,-25){spatial contour function $V(t)$}
\end{picture}

\vskip 10mm
Figure 2
\end{center}

We denote
by
$\Pi(d\t)$ the law of the Galton-Watson tree with offspring distribution
$\mu$, which is  a probability measure on ${\bf T}$. If $\t\in{\bf T}$ and $v_0\in\t$, let
$$\t^{[v_0]}:=\{v\in\u:v_0v\in\t\}$$
denote the subtree of $\t$ originating from $v_0$. 
The probability measure $\Pi$ is characterized by the following two properties \cite{Neveu}:
\begin{description}
\item{(i)} $\Pi(N_\varnothing=j)=\mu(j)$ for every $j\geq 0$.
\item{(ii)} Under the conditional measure $\Pi(\cdot \mid N_\varnothing=j)$, the
subtrees
$\t^{[1]},\t^{[2]},\ldots,\t^{[j]}$ are independent and distributed according to $\Pi$.
\end{description}

The probability measure $\P_x(d\t dU)$ on the space $\Omega$
is then defined by
$$\P_x(d\t dU)=\Pi(d\t)\,R_x(\t,dU)$$
where, for every $\t\in{\bf T}$, the probability measure
$R_x(\t,dU)$ is characterized as follows. Let $\e_\t$ denote the set of
all edges of $\t$ and let $(Y_e,e\in\e_\t)$ be i.i.d. random variables
with distribution $\gamma$. For every $v\in\t$, set
$$X_v=x+\sum_{e\in[\varnothing,v]} Y_e\;,$$
where the notation $e\in[\varnothing,v]$ means that the edge $e$
belongs to the ancestral line of $v$. Then 
$R_x(\t,dU)$ is the distribution of $(X_v,v\in\t)$.

 Let
us a recall a well-known formula for the distribution of 
$|\t|$ under $\Pi$ (see e.g. Section 5.2 of \cite{Pit}). Let $(S_n)_{n\geq 0}$ be a random walk on
$\Z$ with jump distribution $\nu(k)=\mu(k+1)$, $k=-1,0,1,2,\ldots$, 
started from the origin and defined
under the probability measure $P$.
Set $\tau:=\inf\{n\geq 0:S_n=-1\}$. Then, for every integer 
$n\geq 1$,
$$\Pi(|\t|=n)=P(\tau=n)={1\over n}\,P(S_n=-1).$$
The aperiodicity of $\mu$ now implies that the latter quantity is
positive for every $n$ sufficiently large, so that we can define
\ba
&&\Pi^n(d\t)=\Pi(d\t\mid |\t|=n+1)\\
&&\P^n_x(d\t dU)=\P_x(d\t dU\mid |\t|=n+1).
\ea
Later, each time we will consider the probability measures 
$\Pi^n$ or $\P^n_x$, it will be implicit that $n$ is large enough so that
this definition makes sense.

The following result is a special case of Theorem 2 in \cite{JM}. Recall
that $(\eg,Z^0)$ denotes the (one-dimensional) Brownian snake with initial point $0$,
as defined in Section 1.

\begin{theorem}
\label{invar0}
Assume that {\rm (\ref{4mom})} holds. Then the law under $\P^n_0$ of
$$\left( \Big({\sigma\over 2} \;{C(2nt)\over n^{1/2}}\Big)_{0\leq t\leq
1},
\Big({1\over \rho} \Big({\sigma\over 2}\Big)^{1/2}
{V(2nt)\over n^{1/4}}\Big)_{0\leq t\leq 1}\right)$$
converges as $n\to \infty$ to the law 
of $(\eg,Z^0)$. 
The convergence holds in the sense of weak convergence 
of probability measures on $C([0,1],\R)^2$.
\end{theorem}

We aim at proving a conditional version of Theorem \ref{invar0}.
If $(\t,U)$ is a spatial tree, we set
$$\un U=\inf\{U_v:v\in\t,\,v\not =\varnothing\}$$
with the convention that $\un U=+\infty$ if $\t=\{\varnothing\}$.
For every $x\geq 0$, we then define
$$\ov \P^n_x(\cdot):=\P^n_x(\cdot\mid \un U>0).$$

Recall from Section 1 the notation $(\ov \eg^0,\ov Z^0)$ for the conditioned Brownian snake. 
In the notation of
\cite{LGW}, the distribution of $(\ov \eg^0,\ov Z^0)$ is the law of the pair $(\zeta,\wh W)$
under $\ov\N^{(1)}_0$.

\begin{theorem}
\label{invar}
Assume that {\rm (\ref{4mom})} holds and let $x\geq 0$. Then the law under $\ov\P^n_x$
of
$$\left( \Big({\sigma\over 2} \;{C(2nt)\over n^{1/2}}\Big)_{0\leq t\leq
1},
\Big({1\over \rho} \Big({\sigma\over 2}\Big)^{1/2}
{V(2nt)\over n^{1/4}}\Big)_{0\leq t\leq 1}\right)$$
converges as $n\to \infty$ to the law 
of $(\ov \eg^0,\ov Z^0)$. 
The convergence holds in the sense of weak convergence 
of probability measures on $C([0,1],\R)^2$.
\end{theorem}

\noindent{\bf Remark.} A trivial translation argument shows that 
$\P^n_0$ in Theorem \ref{invar0} could be replaced by 
$\P^n_x$ for any $x\in \R$. In the setting of Theorem \ref{invar} however,
no obvious argument can be used to reduce the proof to one particular
value of $x$.

\section{Re-rooting spatial trees}

Recall that
$$\u=\bigcup_{n=0}^\infty \N^n$$
denotes the set of all possible vertices in our discrete trees, 
and that $\u^*=\u\backslash\{\varnothing\}$.

Let $v_0\in\u^*$ and let $\t\in{\bf T}$ such that $v_0\in\t$. Let
$k=k(v_0,\t)$ be the time of the first visit of $v_0$ in the evolution of the
contour of $\t$. Also let $\ell=\ell(v_0,\t)$ be the time of the
last visit of $v_0$. Note that $\ell\geq k$ and $\ell=k$ iff
$v_0$ is a leaf of $\t$. To simplify notation, we set
$\zeta(\t)=2(|\t|-1)$.

For every $t\in[0,\zeta(\t)-(\ell-k)]$, set
$$\wh C^{(v_0)}(t)
=C(k)+C(\llbracket k-t\rrbracket)
-2 \inf_{\llbracket k-t\rrbracket\wedge k\leq t\leq 
\llbracket k-t\rrbracket \vee k} C(r)$$
where $C(\cdot)$ is as above the contour function of $\t$, and
$\llbracket k-t\rrbracket$ stands for the unique element of
$[0,\zeta(\t))$ such that $\llbracket k-t\rrbracket-(k-t)=0$
or $\zeta(\t)$. We also set $\wh C^{(v_0)}(t)=0$ if $t>\zeta(\t)-(\ell-k)$.

Then it is easy to verify that there exists a unique plane tree
$\wh \t^{(v_0)}\in{\bf T}$ whose contour function is $\wh C^{(v_0)}$. Informally,
$\wh \t^{(v_0)}$ is obtained by removing all vertices that are descendants
of $v_0$ and then re-rooting the resulting tree at $v_0$ (we should also
specify the ordering of children in the re-rooted tree, but we omit details
in this informal description). See Fig.3 for an example.

\begin{center}
\unitlength=1.1pt
\begin{picture}(200,125)

\thicklines \put(10,0){\line(-1,2){20}}
\thicklines \put(10,0){\line(1,2){20}}
\thicklines \put(30,40){\line(0,1){40}}
\thicklines \put(30,40){\line(-1,2){20}}
\thicklines \put(30,40){\line(1,2){20}}
\thicklines \put(30,40){\line(1,1){40}}
\thicklines \put(30,80){\line(-1,2){20}}
\thicklines \put(30,80){\line(1,2){20}}
\thicklines \put(10,0){\line(0,1){40}}
\thicklines \put(30,80){\line(0,1){40}}

\put(20,76){$v_0$}
\put(-18,36){$\ov v$}

\thicklines \put(170,0){\line(0,1){40}}
\thicklines \put(170,40){\line(0,1){40}}
\thicklines \put(170,40){\line(-1,2){20}}
\thicklines \put(170,40){\line(1,2){20}}
\thicklines \put(170,40){\line(1,1){40}}
\thicklines \put(170,80){\line(1,2){20}}
\thicklines \put(170,80){\line(-1,2){20}}

\put(160,76){$\wh v_0$}
\put(194,116){$v$}

\put(0,-15){tree $\t$}
\put(140,-15){re-rooted tree $\wh\t^{(v_0)}$}

\end{picture}

\vspace{5mm}
Figure 3
\end{center}

We note that $|\wh \t^{(v_0)}|=|\t|$ iff $v_0$ is a leaf of $\t$.
Also, if $v_0=(j_1,j_2,\ldots,j_p)$, then $\wh v_0:=(1,j_p,j_{p-1},\ldots,j_2)$ 
automatically belongs to $\wh \t^{(v_0)}$. Indeed, $\wh v_0$ is the vertex of 
the re-rooted tree corresponding to the root of the initial tree. 
In Fig.3, $v_0=(3,2)$ and $\wh v_0=(1,2)$.

By definition, $\varnothing$ has exactly one child in the re-rooted
tree $\wh \t^{(v_0)}$. We define a new probability measure $Q$
on $\bf T$ by setting
$$Q(d\t)=\Pi(d\t\mid N_\varnothing =1),$$
where $N_\varnothing$ is the number of children of $\varnothing$
in $\t$. The conditioning a priori makes sense only if $\mu(1)>0$. However, even when
$\mu(1)=0$, there is an obvious way of defining $Q$. 

If $\t\in{\bf T}$ and $w_0\in \t$, we also denote by $\t^{(w_0)}$
the new tree obtained from $\t$ by removing those vertices which 
are descendants of $w_0$ not equal to $w_0$.

\begin{lemma}
\label{re-root-1}
Let $v_0\in\u^*$ of the form $v_0=(1,j_2,\ldots,j_p)$
for some $p\geq 1$, $j_2,\ldots,j_p\in\N$. Assume that 
$Q(v_0\in\t)>0$. Then the law under $Q(\cdot\mid v_0\in\t)$ 
of the re-rooted tree $\wh \t^{(v_0)}$ coincides with the law
under $Q(\cdot \mid \wh v_0\in\t)$ of the tree $\t^{(\hat v_0)}$.
\end{lemma}

The proof is an elementary application of properties of 
Galton-Watson trees. We leave details to the reader.

We shall need a spatial version of Lemma \ref{re-root-1}.
Let $\Q$ be the probability measure on $\Omega$
defined by
$$\Q(d\t dU)=\P_0(d\t dU\mid N_\varnothing =1) =Q(d\t)\,R_0(\t,dU).$$
If $(\t,U)\in\Omega$ and $v_0\in\t\backslash\{\varnothing\}$, the re-rooted
spatial tree
$(\wh \t^{(v_0)}, \wh U^{(v_0)})$ is defined as follows: For every
vertex $v$ of $\wh \t^{(v_0)}$, $\wh U^{(v_0)}_v=U_{\ov v}-U_{v_0}$, if
$\ov v$ is the vertex of the initial tree $\t$ corresponding 
to $v$ in $\wh \t^{(v_0)}$ (see Fig.3 for an example). Alternatively, we may say that the 
spatial contour function $\wh V^{(v_0)}$ of $(\wh \t^{(v_0)}, \wh U^{(v_0)})$
is determined by
$$\wh V^{(v_0)}(t)=V(\llbracket k-t\rrbracket)-V(k)$$
for $0\leq t\leq \zeta(\wh \t^{(v_0)})=\zeta(\t)-(\ell-k)$. (Here 
$k=k(v_0,\t)$ and $\ell=\ell(v_0,\t)$ are as in the beginning of the
section.)

\begin{lemma}
\label{re-root-2}
Let $v_0\in\u^*$ of the form $v_0=(1,j_2,\ldots,j_p)$
for some $p\geq 1$, $j_2,\ldots,j_p\in\N$. Assume that 
$Q(v_0\in\t)>0$. Then the law under $\Q(\cdot\mid v_0\in\t)$ 
of the re-rooted tree $(\wh \t^{(v_0)}, \wh U^{(v_0)})$ coincides with the law
under $\Q(\cdot \mid \wh v_0\in\t)$ of the spatial tree $(\t^{(\hat v_0)},
U^{(\hat v_0)})$, where $U^{(\hat v_0)}$ denotes the restriction of 
$U$ to $\t^{(\hat v_0)}$.
\end{lemma}

Lemma \ref{re-root-2} is a simple consequence of Lemma \ref{re-root-1}
and our definitions. Note that we use the symmetry of the spatial distribution 
$\gamma$.

If $(\t,U)$ is a spatial tree, we denote by $\Delta=\Delta(\t,U)$ the set
of all vertices with minimal spatial position:
$$\Delta=\{v\in\t:U_v=\min_{w\in\t} U_w\}.$$
We also denote by $v_m$ the first element of $\Delta$ in lexicographical order.
Finally, we use the notation $\partial T$ for the set of all leaves of $\t$.

\begin{lemma}
\label{key-re-root}
For any nonnegative measurable functional $F$ on $\Omega$,
$$\Q\Big(F(\wh \t^{(v_m)}, \wh U^{(v_m)}){\bf 1}_{\{|\Delta|=1\,,\,v_m\in\partial
\t\}}\Big)=\Q\Big(F(\t,U)\,|\partial\t|\,{\bf 1}_{\{\un U>0\}}\Big).$$
\end{lemma}

Loosely speaking, this lemma says that the spatial tree re-rooted at the
(first) vertex with minimal spatial position is closely related to the initial
tree conditioned to have positive spatial positions. Compare with Theorem
1.2 in \cite{LGW}.

\proof Let $v_0\in\u^*$ such that $Q(v_0\in\t)>0$. Then
\be
\label{reroot1}
\Q\Big(F(\wh \t^{(v_0)}, \wh U^{(v_0)}){\bf 1}_{\{\Delta=\{v_0\}\}}\Big)
=\Q\Big(F(\wh \t^{(v_0)}, \wh U^{(v_0)}){\bf 1}_{\{
v_0\in\t\,;\,U_v>U_{v_0},\,\forall v\in \t\backslash \{v_0\}\}}\Big).
\ee
Recall from Section 2 the notation $\t^{[v_0]}:=\{v\in \u:v_0v\in \t\}$
for the subtree of $\t$ originating from $v_0$. For every
$v\in \t^{[v_0]}$, put
$$U^{[v_0]}_v:=U_{v_0v}-U_{v_0}.$$
Plainly, under the probability measure $\Q(\cdot \mid v_0\in\t)$,
the spatial tree $(\t^{[v_0]},U^{[v_0]})$ is independent of
$(\t^{(v_0)},U^{(v_0)})$ and has distribution $\P_0$. Since 
$(\wh\t^{(v_0)},\wh U^{(v_0)})$ is 
by construction a function of $(\t^{(v_0)},U^{(v_0)})$,
we can rewrite formula (\ref{reroot1}) as follows:
\be
\label{reroot2}
\Q\Big(F(\wh \t^{(v_0)}, \wh U^{(v_0)}){\bf 1}_{\{\Delta=\{v_0\}\}}\Big)
=\Q\Big(F(\wh \t^{(v_0)}, \wh U^{(v_0)}){\bf 1}_{\{
v_0\in\t\,;\,U_v>U_{v_0},\forall v\in \t^{(v_0)}\backslash \{v_0\}\}}\Big)\,
\P_0(\un U>0).
\ee
Using Lemma \ref{re-root-2}, we then get
\ba
\Q\Big(F(\wh \t^{(v_0)}, \wh U^{(v_0)}){\bf 1}_{\{
v_0\in\t\,;\,U_v>U_{v_0},\forall v\in \t^{(v_0)}\backslash \{v_0\}\}}\Big)
\!\!&=&\!\!\Q\Big(F(\wh \t^{(v_0)}, \wh U^{(v_0)}){\bf 1}_{\{
v_0\in\t\,;\,\hat U^{(v_0)}_v>0,\,\forall v\in \hat\t^{(v_0)}\backslash
\{\varnothing\}\}}\Big)\\
\!\!&=&\!\!\Q\Big(F(\t^{(\hat v_0)}, U^{(\hat v_0)}){\bf 1}_{\{
\hat v_0\in\t\,;\, U^{(\hat v_0)}_v>0,\,\forall v\in \t^{(\hat v_0)}\backslash
\{\varnothing\}\}}\Big).
\ea
Now notice that $(\t^{(v)},U^{(v)})=(\t,U)$ if $v\in\partial \t$.
Moreover, the event $\{v\in\partial \t\}$
is independent of $(\t^{(v)},U^{(v)})$ under $\Q(\cdot \mid v\in \t)$, and
has probability $\mu(0)$. Combining these observations, we get
$$
\Q\Big(F(\wh \t^{(v_0)}, \wh U^{(v_0)}){\bf 1}_{\{
v_0\in\t\,;\,U_v>U_{v_0},\forall v\in \t^{(v_0)}\backslash \{v_0\}\}}\Big)
={1\over \mu(0)}\,
\Q\Big(F(\t, U){\bf 1}_{\{
\hat v_0\in\partial\t\,;\, U_v>0,\,\forall v\in \t\backslash
\{\varnothing\}\}}\Big).$$
Using (\ref{reroot2}) and the preceding equalities, we get
\be
\label{reroot3}
\Q\Big(F(\wh \t^{(v_0)}, \wh U^{(v_0)}){\bf 1}_{\{\Delta=\{v_0\}\}}\Big)
={\P_0(\un U>0)\over \mu(0)}
\;\Q\Big(F(\t, U){\bf 1}_{\{
\hat v_0\in\partial\t\,;\, U_v>0,\,\forall v\in \t\backslash
\{\varnothing\}\}}\Big).
\ee
From the property stated just before (\ref{reroot2}), we easily see that
under the conditional measure $\Q(\cdot \mid \Delta=\{v_0\})$,
the spatial tree $(\t^{[v_0]},U^{[v_0]})$ is independent of 
$(\t^{(v_0)},U^{(v_0)})$ and has distribution $\P_0(\cdot\mid \un U>0)$.
Hence,
\be
\label{reroot4}
\Q\Big(F(\wh \t^{(v_0)}, \wh U^{(v_0)}){\bf 1}_{\{\Delta=\{v_0\}\}}\Big)
={1\over \P_0(\t=\{\varnothing\}\mid \un U>0)}
\Q\Big(F(\wh \t^{(v_0)}, \wh U^{(v_0)}){\bf 1}_{\{
v_0\in\partial \t,\,\Delta=\{v_0\}\}}\Big).
\ee
Since
$$\P_0(\t=\{\varnothing\}\mid \un U>0)={\mu(0)\over \P_0(\un U >0)},$$
(\ref{reroot3}) and (\ref{reroot4}) give
\be
\label{reroot5}
\Q\Big(F(\wh \t^{(v_0)}, \wh U^{(v_0)}){\bf 1}_{\{
v_0\in\partial \t,\,\Delta=\{v_0\}\}}\Big)
=\Q\Big(F(\t, U){\bf 1}_{\{
\hat v_0\in\partial\t,\, \un U>0\}}\Big).
\ee
Summing (\ref{reroot5}) over all possible choices of $v_0$ leads to the desired
result. \cq

\smallskip
We shall need a variant of Lemma \ref{key-re-root}.

\begin{lemma}
\label{key-re-root2}
For any nonnegative measurable functional $F$ on $\Omega$,
$$\Q\Big(\sum_{v_0\in \Delta\cap\partial\t} F(\wh \t^{(v_0)}, \wh
U^{(v_0)})\Big)=\Q\Big(F(\t,U)\,|\partial\t|\,{\bf 1}_{\{\un U\geq 0\}}\Big).$$
\end{lemma}

\proof By arguing as in the proof of (\ref{reroot3}), we get for every 
$v_0\in\u^*$,
\ba
\Q\Big(F(\wh \t^{(v_0)}, \wh U^{(v_0)}){\bf 1}_{\{v_0\in\Delta\}}\Big)
\!\!&=&\!\!
\Q\Big(F(\wh \t^{(v_0)}, \wh U^{(v_0)}){\bf 1}_{\{
v_0\in\t\,;\,U_v\geq U_{v_0},\,\forall v\in \t^{(v_0)}\}}\Big)\,\P_0(\un U\geq 0)\\
\!\!&=&\!\!
\Q\Big(F(\t^{(\hat v_0)}, U^{(\hat v_0)}){\bf 1}_{\{
\hat v_0\in\t\,;\, U^{(\hat v_0)}_v\geq 0,\,\forall v\in \t^{(\hat v_0)}\}}\Big)\,
\P_0(\un U\geq 0)\\
\!\!&=&\!\!{\P_0(\un U\geq 0)\over \mu(0)}\;
\Q\Big(F(\t, U){\bf 1}_{\{
\hat v_0\in\partial\t\,;\, U_v\geq 0,\,\forall v\in \t\}}\Big).
\ea
On the other hand, analogously to (\ref{reroot4}),
$$\Q\Big(F(\wh \t^{(v_0)}, \wh U^{(v_0)}){\bf 1}_{\{v_0\in\Delta\}}\Big)
={1\over \P_0(\t=\{\varnothing\}\mid \un U\geq 0)}
\Q\Big(F(\wh \t^{(v_0)}, \wh U^{(v_0)})\,{\bf 1}_{\{
v_0\in\Delta\cap \partial \t\}}\Big).
$$
The lemma follows by combining the previous two identities and then summing
over all choices of $v_0\in\u^*$. \cq

\smallskip
\noindent{\bf Remark.} If $\gamma$ has no atoms we have $|\Delta|=1$,
$\Q$ a.s., and Lemmas \ref{key-re-root} and \ref{key-re-root2} reduce to the
same identity. In our applications, we shall be concerned with the case when 
$\gamma$ does have atoms, and is even supported on a finite subset
of $\Z$.

\section{Estimates for the probability of staying on the positive side}

Our goal in this section is to derive upper and lower bounds for the
probability $\P^n_x(\un U>0)$ when $n\to\infty$. Our main tools
will be Lemmas \ref{key-re-root} and \ref{key-re-root2}. We also need a
preliminary
estimate concerning the cardinality $|\partial \t|$ of the set of leaves.

\begin{lemma}
\label{leaves}
There exists a constant $\alpha_0>0$ such that, for every $n$ sufficiently large,
$$\Q(||\partial\t|-\mu(0)n|>n^{3/4},|\t|=n+1)\leq e^{-n^{\alpha_0}}.$$
\end{lemma}

\proof We first recall some classical facts about relations between random
walks and Galton-Watson trees. As in Section 2 above, 
let $S_n=X_1+\cdots+X_n$ be a random walk on $\Z$
with jump distribution $\nu(k)=\mu(k+1)$, $k=-1,0,1,2,\ldots$, 
started from the origin and defined
under the probability measure $P$.
Set $\tau:=\inf\{n\geq 0:S_n=-1\}$, and for every integer $n\geq 1$,
$$M_n=|\{k\in\{1,2,\ldots,n\}:X_k=-1\}|.$$
Then the law of the pair $(|\t|,|\partial \t|)$ under $Q$
coincides with that of the pair $(1+\tau,M_\tau)$ under $P$. For a proof, see e.g.
the discussion in Section 2 of \cite{LGLJ1}, or Section 5.2 of \cite{Pit}.

It follows that
$$\Q(||\partial\t|-\mu(0)n|>n^{3/4},|\t|=n+1)
=P(|M_\tau-\mu(0)n|>n^{3/4},\tau=n)
\leq\! P(|M_n-\mu(0)n|>n^{3/4}).$$
The estimate of the lemma now follows from standard moderate deviations
estimates for sums of independent Bernoulli variables. \cq

\smallskip
We set $\Q^n(\cdot)=\Q(\cdot\mid|\t|=n+1)$, which makes sense
for every $n$ sufficiently large.

\begin{proposition}
\label{positive}
Let $K>0$. There exist positive constants $c_1,c_2=c_2(K),\wt c_1,\wt c_2$
such that, for every $x\in[0,K]$ and every $n$ sufficiently large,
\ba
&&{\wt c_1\over n}\leq \Q^n(\un U>0) \leq {\wt c_2\over n}\,,\\
\noalign{\smallskip}
&&{c_1\over n}\leq \P^n_x(\un U>0) \leq {c_2\over n}\,.
\ea
\end{proposition}

\noindent{\bf Remark.} By an obvious comparison argument, $c_1$
can be chosen independently of $K$. 

\proof  We first bound $\Q^n(\un U>0)$. We apply Lemma \ref{key-re-root} with 
$$F(\t,U)={\bf 1}_{\{|\t|=n+1\}},$$
noting that $|\wh\t^{(v)}|=|\t|$ if $v\in\partial \t$. It follows that
$$\Q(|\partial \t|\,{\bf 1}_{\{|\t|=n+1\,,\,\un U>0\}})\leq 
\Q(|\t|=n+1).$$

On the other hand, Lemma \ref{leaves} shows that, for $n$ 
sufficiently large,
$$\Q(|\partial \t|\,{\bf 1}_{\{|\t|=n+1\,,\,\un U>0\}})
\geq (\mu(0)n-n^{3/4})\,\Q(|\t|=n+1,\un U>0)-n\,e^{-n^{\alpha_0}}.$$
By combining this with the preceding bound, we get
\be
\label{posi1}
\Q^n(\un U>0)\leq {1\over \mu(0)n-n^{3/4}}
+{n\,e^{-n^{\alpha_0}}\over (\mu(0)n-n^{3/4})\Q(|\t|=n+1)}.
\ee
With the notation of the proof of Lemma \ref{leaves}, we have also
\be
\label{posi2}
\Q(|\t|=n+1)=P(\tau=n)={1\over n}\,P(S_n=-1)\build{\sim}_{n\to\infty}^{}
c_0\,n^{-3/2}
\ee
where $c_0$ is a positive constant, and the last estimate follows from a
standard local limit theorem. We thus deduce from (\ref{posi1}) that
\be
\label{posi3}
\limsup_{n\to \infty} n\,\Q^n(\un U>0)\leq {1\over \mu(0)},
\ee
yielding the desired upper bound for $\Q^n(\un U>0)$.

Let us now discuss a lower bound for $\Q^n(\un U\geq 0)$. Applying Lemma
\ref{key-re-root2} with the same function $F$, we get
\be
\label{posi4}
\Q(|\partial \t|\,{\bf 1}_{\{|\t|=n+1\,,\,\un U\geq 0\}})
=\Q(|\Delta\cap \partial \t|\,{\bf 1}_{\{|\t|=n+1\}})
\geq \Q(|\Delta\cap \partial \t|\geq 1\,,\,|\t|=n+1),
\ee
and we now need to show that the latter quantity is 
bounded below by $c\,\Q(|\t|=n+1)$ for some
positive constant $c$. To this end, we state another lemma. Recall
that $v_m$ is the first vertex in $\Delta$ for the lexicographical order 
of vertices.

\begin{lemma}
\label{minleaves}
There exists a constant $\ov c_1>0$ such that, for every $n$ sufficiently large,
$$\Q^n(v_m\in\partial \t)\geq \ov c_1.$$
\end{lemma}

We postpone the proof of Lemma \ref{minleaves} and complete that
of Proposition \ref{positive}. As a consequence of Lemma \ref{minleaves},
we have
$$\Q(|\Delta\cap \partial \t|\geq 1\,,\,|\t|=n+1)
\geq \ov c_1\,\Q(|\t|=n+1).$$
Using now (\ref{posi4}), we get
$$n\,\Q(|\t|=n+1\,,\,\un U\geq 0)\geq \ov c_1\,\Q(|\t|=n+1),$$
and so
\be
\label{posi5}
\Q^n(\un U\geq 0)\geq {\ov c_1\over n}.
\ee

The bounds involving $\P^n_x(\un U>0)$ are easily derived 
from (\ref{posi3}) and (\ref{posi5}). Consider first the lower bound. Clearly,
it is enough to take $x=0$. Then fix $k\geq 1$ with $\mu(k)>0$.
By arguing on an appropriate event contained in $\{N_\varnothing=k,
N_1=k\}$, we get
$$\P_0(\un U>0\,,\,|\t|=n+1)\geq
\mu(k)^2\,\gamma((0,\infty))^{2k-1}\,\mu(0)^{2(k-1)}\,\Q(\un U\geq 0,
|\t|=n+2-2k).$$
Since $\P_0(|\t|=n+1)\sim c_0\,n^{-3/2}\sim \Q(|\t|=n+2-2k)$ as $n\to \infty$,
we see that the lower bound for $\P^n_0(\un U>0)$ follows from
(\ref{posi5}).

Consider then the upper bound. It is enough to take $x=K$. Let
$k$ be as above, then choose $y>0$ such that $\gamma([y,\infty))>0$
and let $p\geq 1$ be an integer such that $py\geq K$. Again by arguing
on an appropriate event contained in $\{N_1=k,N_{(1,1)}=k,\ldots,N_{1^{p-1}}=k\}$, we get
\ba &&\Q(\un U>0\,,\,|\t|=n+2+(p-1)k)\\
\noalign{\smallskip}
&&\qquad \geq \mu(k)^{p-1}\gamma([y,\infty))^p\,\mu(0)^{(p-1)(k-1)}\,
\gamma((0,\infty))^{(p-1)(k-1)}\,\P_K(\un U>0,|\t|=n+1)
\ea
and thus the upper bound for $\P^n_K(\un U>0)$ follows from
(\ref{posi3}). 

Finally, the bound $\P^n_0(\un U>0)\geq c_1/n$ readily 
implies $\Q^n(\un U>0)\geq \wt c_1/n$ with $\wt c_1=\gamma((0,\infty))\,c_1$. \cq

\smallskip
\noindent{\bf Proof of Lemma \ref{minleaves}:} We first observe that under the
probability measure $\Q(\cdot\mid v_m\not =\varnothing)$, the spatial tree
$(\t^{[v_m]},U^{[v_m]})$ is independent of $(\t^{(v_m)},U^{(v_m)})$
and has distribution $\P_0(\cdot\mid \un U\geq 0)$. Indeed, if we condition
on the value of $v_m$, we see that this statement follows from
basic properties of Galton-Watson trees, similar to those that were used
in the proof of Lemma \ref{key-re-root}. We then get
\begin{eqnarray}
\label{minleaves1}
\Q(|\t|=n+1\,,\,v_m\not =\varnothing)
\!\!&=&\!\!
\sum_{k=0}^{n-1} \Q(|\t^{(v_m)}|=n-k+1\,,\,|\t^{[v_m]}|=k+1\,,\,v_m\not
=\varnothing)\nonumber\\
\!\!&=&\!\!\sum_{k=0}^{n-1} \Q(|\t^{(v_m)}|=n-k+1\,,\,v_m\not
=\varnothing)\,\eta(k+1)
\end{eqnarray}
where $\eta(j)=\P_0(|\t|=j\mid \un U\geq 0)$, for every $j\geq 1$. On the other
hand, using once again the observation of the beginning of the proof, we have
for every integer $\ell\geq 2$,
\begin{eqnarray}
\label{minleaves2}
\Q(|\t|=\ell\,,\,v_m\in\partial \t)
\!\!&=&\!\!\Q(|\t^{(v_m)}|=\ell\,,\,\t^{[v_m]}=\{\varnothing\}\,,\,v_m\not =
\varnothing)\nonumber\\
\!\!&=&\!\! \eta(0)\,\Q(|\t^{(v_m)}|=\ell\,,\,v_m\not =\varnothing).
\end{eqnarray}
Fix an integer $p\geq 2$ such that $\mu(p)>0$. Under the probability
measure $\Q(\cdot\mid N_1=p)$, we can
consider the $p$ trees $\t_1,\ldots,\t_p$
defined as follows. For every $1\leq j\leq p$, $\t_j$ consists of the
root $\varnothing$ and of the vertices of the type $1v$
such that $1jv\in\t$. Then, under $\Q(\cdot \mid N_1=p)$, 
$\t_1,\ldots,\t_p$ are independent and distributed according to $Q$.
Let $k$ be an integer with $p-1\leq k\leq n-2$. We consider the
event where $|\t_1|=n-k$, $|\t_2|=k-p+3$
and $|\t_j|=2$ for other values of $j$. By arguing on this event
and imposing appropriate conditions on the spatial displacements
(using in particular the fact that $\gamma((-\infty,0])=\gamma([0,\infty))
\geq 1/2$), we get the bound
\ba
&&\Q(|\t|=n+1\,,\,v_m\in\partial \t)\\
&&\qquad\geq{\mu(p)\over 2} \sum_{k=p-1}^{n-2} 
\Q(|\t|=n-k\,,\,v_m\in\partial\t)\,
\Q(|\t|=k-p+3\,,\,\un U\geq 0)\,\Big({\mu(0)\over 2}\Big)^{p-2}\\
&&\qquad =c_{(p)}\,\sum_{k=p-1}^{n-2} 
\Q(|\t|=n-k\,,\,v_m\in\partial\t)\,
\Q(|\t|=k-p+3\,,\,\un U\geq 0)
\ea
where $c_{(p)}>0$ is a constant depending on $p$ and $\mu$. 

Similarly, by requiring the spatial displacement along the first edge to be
nonnegative, we get for $k\geq p-1$,
$$\Q(|\t|=k-p+3\,,\,\un U\geq 0)
\geq {1\over 2}\,\P_0(|\t|=k-p+2\,,\,\un U\geq 0)
={1\over 2}\,\P_0(\un U\geq 0)\,\eta(k-p+2).$$
It follows that
\ba
\Q(|\t|=n+1\,,\,v_m\in\partial \t)
\!\!&\geq&\!\! c'_{(p)}\,\sum_{k=p-1}^{n-2} 
\Q(|\t|=n-k\,,\,v_m\in\partial\t)\,\eta(k-p+2)\\
\!\!&=&\!\! c'_{(p)}\,\sum_{j=0}^{n-p-1} 
\Q(|\t|=(n-p)-j+1\,,\,v_m\in\partial\t)\,\eta(j+1).
\ea
Using (\ref{minleaves2}) we arrive at
\be
\label{minleaves3}
\Q(|\t|=n+1\,,\,v_m\in\partial \t)
\geq c''_{(p)}\,\sum_{j=0}^{(n-p)-1} 
\Q(|\t^{(v_m)}|=(n-p)-j+1\,,\,v_m\not =\varnothing)\,\eta(j+1).
\ee
We compare the last bound with (\ref{minleaves1}) written with
$n$ replaced by $n-p$. It follows that
\be
\label{minleaves4}
\Q(|\t|=n+1\,,\,v_m\in\partial \t)\geq c''_{(p)}\,
\Q(|\t|=n-p+1\,,\,v_m\not =\varnothing).
\ee
Since $\Q(|\t|=n+1\,,\,v_m\not =\varnothing)\sim \Q(|\t|=n+1)\sim c_0n^{-3/2}$
as $n\to \infty$, Lemma \ref{minleaves} follows from (\ref{minleaves4}). \cq

\smallskip
\noindent{\bf Remark} In view of Proposition \ref{positive}, one expects the
existence of a positive constant $c_\infty$ such that
$$\lim_{n\to \infty} n\,\Q^n(\un U>0)=c_\infty.$$
When $\gamma$ has no atoms, and under some additional conditions on 
the offspring distribution $\mu$, this can indeed be proved from the
lemmas of Section 3, with $c_\infty=\P_0(\un U>0)^{-1}$. Similar asymptotics
for $\P^n_0(\un U>0)$ then follow easily. Since we do not
need these more precise estimates, we will not address this problem here.

\section{A spatial Markov property}

In this section, we briefly discuss a Markov property for our branching trees,
which will be used in the proof of our main result. Arguments are elementary
and so we omit most details.

We fix $a>0$. If $(\t,U)$ is a spatial tree and $v\in\t$, we say that $v$ is
an exit vertex from $(-\infty,a)$ if $U_v\geq a$ and $U_{v'}<a$
for every ancestor $v'$ of $v$ distinct from $v$. Denote by $v_1,v_2,\ldots,v_M$ 
the exit vertices from $(-\infty,a)$ listed in lexicographical order.
If $i,j\in\{1,2,\ldots,M\}$ and $i\not =j$, then $v_i$ cannot be an
ancestor of $v_j$. 

For $v\in\t$, and for every $v'\in\t^{[v]}$, we set
$$\ov U^{[v]}_{v'}=U_{vv'}$$
(compare with the definition of $U^{[v]}$). Finally, we denote by $\t^a$
the subtree of $\t$ consisting of those vertices which are not
strict descendants of $v_1,\ldots,v_M$. In particular, $v_1,\ldots,v_M\in\t^a$.
We also denote by $U^a$ the restriction of $U$ to $\t^a$. Informally,
$(\t^a,U^a)$ corresponds to the tree $(\t,U)$ ``truncated at the first exit
time'' from $(-\infty,a)$.

\begin{proposition}
\label{spatial0}
Let $x\in[0,a)$ and $p\geq 1$. Under the probability measure $\P_x(\cdot\mid
M=p)$, conditionally on $(\t^a,U^a)$, the spatial trees $(\t^{[v_1]},\ov
U^{[v_1]}),\ldots,(\t^{[v_p]},\ov U^{[v_p]})$ are independent and distributed
respectively according to $\P_{U_{v_1}},\ldots,\P_{U_{v_p}}$.
\end{proposition}

The proof of Proposition \ref{spatial0} is an easy application of
properties of Galton-Watson trees. We leave details to the reader.
See e.g. \cite{Ch} for closely related statements in a slightly
different setting.

Conditional versions of Proposition \ref{spatial0} are derived in
a straightforward way. Firstly, this statement remains valid if
$\P_x$ is replaced by
$$\ov \P_x(\cdot):=\P_x(\cdot\mid \un U>0),$$
provided $\P_{U_{v_1}},\ldots,\P_{U_{v_p}}$ in the conclusion are also replaced 
by $\ov\P_{U_{v_1}},\ldots,\ov\P_{U_{v_p}}$.

Then, by conditioning with respect to the sizes of the various trees, we
arrive at the following result. 

\begin{corollary}
\label{spatial}
Let $x\in[0,a)$ and $p\in\{1,\ldots,n\}$. Let $n_1,\ldots,n_p$ be positive
integers such that $n_1+\cdots+n_p\leq n$. Assume that
$$\ov\P^n_x(M=p,|\t^{[v_1]}|=n_1,\ldots,|\t^{[v_p]}|=n_p)>0.$$
Then, under the probability measure $\ov\P^n_x(\cdot\mid M=p
,|\t^{[v_1]}|=n_1,\ldots,|\t^{[v_p]}|=n_p)$, 
conditionally on $(\t^a,U^a)$, the spatial trees $(\t^{[v_1]},\ov
U^{[v_1]}),\ldots,(\t^{[v_p]},\ov U^{[v_p]})$ are independent and distributed
respectively according to $\ov\P^{n_1}_{U_{v_1}},\ldots,\ov\P^{n_p}_{U_{v_p}}$.
\end{corollary}

\section{Asymptotic properties of conditioned trees}

From now on we assume that (\ref{4mom}) holds.

In view of our main result Theorem \ref{invar}, it is convenient 
to introduce a specific notation for rescaled processes. For every
integer $n\geq 1$ and every $t\in[0,1]$, we set
\ba
&&C^{(n)}(t)={\sigma\over 2} \;{C(2nt)\over n^{1/2}}\,,\\
&&V^{(n)}(t)={1\over \rho}
\Big({\sigma\over 2}\Big)^{1/2} {V(2nt)\over n^{1/4}}\,.
\ea
Before proceeding to the proof of Theorem \ref{invar}, we need to get some
information about asymptotic properties of the pair $(C^{(n)},V^{(n)})$
under $\ov\P^n_x$. We will consider the conditioned measure
$$\ov\Q^n:=\Q^n(\cdot\mid \un U>0).$$

\begin{proposition}
\label{keytech}
For every $b>0$ and $\varepsilon\in(0,1/10)$, 
we can find
$\delta,\alpha\in(0,\varepsilon)$ such that, for all $n$ sufficiently large,
$$\ov\Q^n\Big(\inf_{t\in[\delta/2,1-\delta/2]}V^{(n)}(t)\geq 2\alpha\;,\;
\sup_{t\in[0,4\delta]\cup[1-4\delta,1]}(C^{(n)}(t)+V^{(n)}(t))\leq
\varepsilon/2\Big)
\geq 1-b.$$
Consequently, if $K>0$, we have also for all $n$ sufficiently large, for every $x\in[0,K]$,
$$\ov\P_x^n\Big(\inf_{t\in[\delta,1-\delta]}V^{(n)}(t)\geq \alpha\;,\;
\sup_{t\in[0,3\delta]\cup[1-3\delta,1]}(C^{(n)}(t)+V^{(n)}(t))\leq
\varepsilon\Big)
\geq 1-c_3b,$$
where the constant $c_3$ only depends on $\mu$, $\gamma$ and $K$.
\end{proposition}

The second part of the proposition will follow from the first one by
arguments similar to those that were used in the proof of Proposition
\ref{positive}
above. 
To prove the first part of the proposition, we will use Theorem \ref{invar0}
together with the following crucial lemma.

\begin{lemma}
\label{reroot-equiv}
Let $F$ be a nonnegative measurable function on $\Omega$
such that $0\leq F\leq 1$. There exist a finite constant $\ov c$, which
does not depend on $F$ nor on $n$, such that
$$\ov\Q^n(F(\t,U))\leq \ov c\,\Q^n(F(\wh\t^{(v_m)},\wh U^{(v_m)}))
+O(n^{5/2}e^{-n^{\alpha_0}}),$$
where $\alpha_0$ is as in Lemma \ref{leaves}, and the estimate
$O(n^{5/2}e^{-n^{\alpha_0}})$ for the remainder holds uniformly in $F$.
\end{lemma}

\proof For every $n$ sufficiently large,
\ba
\Q\Big(F(\t,U)\,{\bf 1}_{\{|\t|=n+1\}}\,{\bf 1}_{\{\un U>0\}}\Big)
\!\!&\leq&\!\!
\Q\Big(F(\t,U)\,{\bf 1}_{\{|\t|=n+1\}}\,{\bf 1}_{\{\un U>0\}}
\,{|\partial \t|\over \mu(0)n-n^{3/4}}\,{\bf 1}_{\{
|\partial \t|\geq \mu(0)n-n^{3/4}\}}\Big)\\
&&\qquad+\Q\Big(F(\t,U)\,{\bf 1}_{\{|\t|=n+1\}}\,{\bf 1}_{\{\un U>0\}}
\,{\bf 1}_{\{
|\partial \t|< \mu(0)n-n^{3/4}\}}\Big)\\
\!\!&\leq&\!\! {2\over \mu(0)n}\,\Q\Big(F(\wh\t^{(v_m)},\wh U^{(v_m)})\,
{\bf 1}_{\{|\t|=n+1\}}\Big)+e^{-n^{\alpha_0}},
\ea
using Lemma \ref{key-re-root} and Lemma \ref{leaves} in the last bound.
Dividing by $\Q(|\t|=n+1)$, we get
$$\Q^n(F(\t,U)\,{\bf 1}_{\{\un U>0\}})
\leq {2\over \mu(0)n}\,\Q^n(F(\wh\t^{(v_m)},\wh U^{(v_m)}))+
O(n^{3/2}e^{-n^{\alpha_0}}).$$
By Proposition \ref{positive}, we have 
$\Q^n(\un U>0)\geq \wt c_1/n$. The lemma now follows from the preceding bound, with
$\ov c=2/(\mu(0)\wt c_1)$. \cq

\smallskip
\noindent{\bf Proof of Proposition \ref{keytech}:} {\it First step.}
We first observe that Theorem \ref{invar0} obviously remains valid 
if $\P^n_0$ is replaced by $\Q^n$. By the Skorokhod representation theorem, we
can find, for every integer $n$ sufficiently large, a pair $(C_n,V_n)$ 
such that the following holds. The processes $C_n$ and $V_n$ are respectively the
contour function and the spatial contour function of a spatial tree $(\t_n,U_n)$
with distribution 
$\Q^n$. Moreover,
\be
\label{Sko}
\left( \Big({\sigma\over 2} \;{C_n(2nt)\over n^{1/2}}\Big)_{0\leq t\leq
1},
\Big({1\over \rho} \Big({\sigma\over 2}\Big)^{1/2}
{V_n(2nt)\over n^{1/4}}\Big)_{0\leq t\leq 1}\right)
\build{\la}_{n\to\infty}^{} (\eg,Z^0),
\ee
uniformly on [0,1], a.s., and the limiting pair $(\eg,Z^0)$ 
is the Brownian snake with initial point $0$, as
defined in the introduction above.

In agreement with the previous notation, write $v_m^n$ for the 
first vertex realizing the minimal spatial position in $\t_n$, and 
$k_n$, respectively $\ell_n$, for the first, resp. the last, time of visit
of $v^n_m$ in the evolution of the contour of $\t_n$.

From Proposition 2.5 in \cite{LGW}, we know that there is a.s. a unique 
$s_*\in(0,1)$ such that
$$Z^0(s_*)=\inf_{0\leq t\leq 1} Z^0(t).$$
The convergence (\ref{Sko}) then implies that
\be 
\label{Skomin}
\lim_{n\to\infty} {k_n\over 2n}=\lim_{n\to\infty} {\ell_n\over 2n}=s_*\ ,\quad
\hbox{a.s.}
\ee

Consider then the re-rooted tree $(\wh\t^{(v^n_m)}_n,\wh U^{(v^n_m)}_n)$.
By construction, its contour function is
$$\wh C^{(v^n_m)}_n(t)=C_n(k_n)+C_n(\llbracket k_n-t\rrbracket_n)
-2\,\inf_{\llbracket k_n-t\rrbracket_n\wedge k_n\leq r\leq \llbracket
k_n-t\rrbracket_n\vee k_n} C_n(r),$$
for $0\leq t\leq 2n-(\ell_n-k_n)$. (Here $\llbracket k_n-t\rrbracket_n$
denotes the unique element of $[0,2n)$ such that $\llbracket k_n-t\rrbracket_n
-(k_n-t)=0$ or $2n$.) The corresponding spatial contour function is
$$\wh V^{(v^n_m)}_n(t)=V_n(\llbracket k_n-t\rrbracket_n)-V_n(k_n).$$
For $2n-(\ell_n-k_n)<t\leq 2n$, we also set $\wh C^{(v^n_m)}_n(t)=\wh V^{(v^n_m)}_n(t)
=0$.
From (\ref{Sko}), (\ref{Skomin}), and the preceding formulas 
for $\wh C^{(v^n_m)}_n(t)$ and $\wh V^{(v^n_m)}_n(t)$, we get
\be
\label{Skorev}
\left( \Big({\sigma\over 2} \;{\wh C^{(v^n_m)}_n(2nt)\over n^{1/2}}\Big)_{0\leq
t\leq 1},
\Big({1\over \rho} \Big({\sigma\over 2}\Big)^{1/2}
{\wh V^{(v^n_m)}_n(2nt)\over n^{1/4}}\Big)_{0\leq t\leq 1}\right)
\build{\la}_{n\to\infty}^{} (\ov\eg^0,\ov Z^0),
\ee
uniformly on [0,1], a.s., where, as in Section 1,
\ba
&&\ov\eg^0(t)=\eg(\{s_*-t\})+\eg(s_*)-2\,\inf_{\{s_*-t\}\wedge s_*\leq r\leq 
\{s_*-t\}\vee s_*}\eg(r)\\
&&\ov Z^0(t)=Z^0(\{s_*-t\})-Z^0(s_*),
\ea
where $\{r\}$ denotes the fractional part of $r$. 

Write ${\bf P}$ for the probability measure under which the processes
$(C_n,V_n)$ and $(\eg,Z^0)$ are defined.
From Lemma \ref{reroot-equiv} applied with a suitable indicator function $F$, we
have for every choice of
$\alpha,\delta,\varepsilon>0$,
\begin{eqnarray}
\label{Sko2}
&&\!\!\!\!\!\!\!\!\limsup_{n\to \infty}
\ov\Q^n\Big(\Big\{\inf_{t\in[\delta/2,1-\delta/2]}V^{(n)}(t)< 2\alpha\Big\}\cup
\Big\{\sup_{t\in[0,4\delta]\cup[1-4\delta,1]}(C^{(n)}(t)+V^{(n)}(t))>
{\varepsilon\over 2}\Big\}\Big)\nonumber\\
&&\!\!\!\!\!\!\!\!\leq \ov c \,\limsup_{n\to\infty}
\Q^n\Big(\!\Big\{\!\inf_{t\in[\delta/2,1-\delta/2]}\!\wh V^{(v_m),(n)}(t)\!<\!
2\alpha\Big\}\cup
\Big\{\!\sup_{t\in[0,4\delta]\cup[1-4\delta,1]}\!(\wh
C^{(v_m),(n)}(t)+\wh V^{(v_m),(n)}(t))\!>\!
{\varepsilon\over 2}\Big\}\!\Big)\nonumber\\
&&\!\!\!\!\!\!\!\!\leq \ov c\,{\bf P}\Big(\Big\{\inf_{t\in[\delta/2,1-\delta/2]}
\ov Z^0(t)\leq 2\alpha\Big\}\cup\Big\{
\sup_{t\in[0,4\delta]\cup[1-4\delta,1]} (\ov\eg^0(t)+\ov Z^0(t))\geq {\varepsilon\over 2}\Big\}
\Big),
\end{eqnarray}
where we used the notation
$$\wh C^{(v_m),(n)}(t)={\sigma\over 2} \;{\wh C^{(v_m)}(2nt)\over n^{1/2}}\quad,
\quad \wh C^{(v_m),(n)}(t)={1\over \rho}
\Big({\sigma\over 2}\Big)^{1/2} {\wh V^{(v_m)}(2nt)\over n^{1/4}},
$$
if $0\leq t\leq (|\t^{(v_m)}|-1)/n$, and $\wh C^{(v_m),(n)}(t)=\wh V^{(v_m),(n)}(t)=0$
if $(|\t^{(v_m)}|-1)/n<t\leq 1$.
In the last inequality above, we used (\ref{Skorev}) together with the fact that
$C_n$ and $V_n$ are respectively the contour function and the
spatial contour function of a spatial tree with distribution 
$\Q^n$.

Recall that $\ov Z^0(t)>0$ for every $t\in(0,1)$, a.s. Hence, if $b>0$ and
$\varepsilon>0$
are given, we can first choose $\delta\in(0,\varepsilon)$ so small that
$$\ov c\,{\bf P}\Big(\sup_{t\in[0,4\delta]\cup[1-4\delta,1]} (\ov\eg^0(t)+\ov
Z^0(t))\geq {\varepsilon\over 2}\Big)< {b\over 2}$$
and then find $\alpha\in(0,\varepsilon)$ small enough so that
$$\ov c\,{\bf P}\Big(\inf_{t\in[\delta/2,1-\delta/2]}
\ov Z^0(t)\leq 2\alpha\Big)<{b\over 2}.$$
The first part of Proposition \ref{keytech} then follows from (\ref{Sko2}).

\noindent{\it Second step.} We now explain how the desired bound under
$\ov\P^n_x$ can be deduced from the one under $\ov\Q$. In a way very similar to
the end of the proof of Proposition \ref{positive}, we first choose an 
integer $\ell\geq 1$ such that $\mu(\ell)>0$. Then let $y>0$ be such that
$\gamma((y,y+1))>0$ and let $p\geq 1$ be the first integer such that $py\geq K$. 
Also set $m=n+(p-1)\ell+1$. Under an appropriate conditioning 
of $\Q$ (requiring in particular that $N_1=\ell, N_{(1,1)}=\ell,\ldots,N_{1^{p-1}}=\ell$), we can
embed a tree with distribution $\P_z$, for some random $z\geq K$, into a tree
distributed according to $\Q$, and we arrive at the bound
\begin{eqnarray}
\label{keytech1}
&&\Q\Big(\{|\t|=m+1\}\cap \{\un U>0\}\cap \Big\{\sup_{t\in[0,A+p]\cup
[2n-A,2m]}(C(t)+V(t))\geq \lambda\Big\}\Big)\nonumber\\
&&\qquad\geq \beta\;\P_x\Big(\{|\t|=n+1\}
\cap \{\un U>0\}\cap \Big\{\sup_{t\in[0,A]\cup
[2n-A,2n]}(C(t)+V(t))\geq \lambda\Big\}\Big),
\end{eqnarray}
where
$$\beta=\mu(\ell)^{p-1}\,\gamma((y,y+1))^p\,
\mu(0)^{(p-1)(\ell-1)}\,\gamma((0,\infty))^{(p-1)(\ell-1)}$$
and $x\in[0,K]$, $A\in(0,2n)$, $\lambda>0$ are arbitrary. Similarly, with the same constant
$\beta$, we have
\begin{eqnarray}
\label{keytech2}
&&\Q\Big(\{|\t|=m+1\}\cap \{\un U>0\}\cap
\Big\{\inf_{t\in[A,2m-A]}V(t)\leq \lambda\Big\}\Big)\nonumber\\
&&\qquad\geq \beta\;\P_x\Big(\{|\t|=n+1\}
\cap \{\un U>0\}\cap \Big\{\inf_{t\in[A,2n-A]}V(t)\leq \lambda-(K+y+1)\Big\}\Big).
\end{eqnarray}
Also recall that the quantities $\Q(\{|\t|=m+1\}\cap \{\un U>0\})$ 
and $\P_x(\{|\t|=n+1\}\cap \{\un U>0\})$ are bounded above and below 
by positive constants times $n^{-5/2}$. Using this last remark, we see that
the second part of Proposition \ref{keytech} follows from the first part,
(\ref{keytech1}) and (\ref{keytech2}). \cq

\medskip
\noindent{\bf Remark.} The limiting process $(\ov\eg^0,\ov Z^0)$ in (\ref{Skorev})
is the same as the one in Theorem \ref{invar}. Therefore, it seems tempting to deduce Theorem
\ref{invar} from (\ref{Skorev}) and the relations between the conditioned spatial tree
and the tree re-rooted at its first minimum (Lemmas \ref{key-re-root}
and \ref{key-re-root2}). This approach would indeed be successful, maybe under 
additional assumptions, in the case when the probability measure $\gamma$
has no atoms, so that the minimal spatial position is attained at a
unique vertex. In our general setting however, we will have to
use a different argument which is explained in the next section.

\section{Proof of the main result}

In this section, we prove Theorem \ref{invar}. We equip $C([0,1],\R)^2$
with the norm $\|(f,g)\|=\|f\|_u\vee\|g\|_u$, where $\|f\|_u$
stands for the uniform norm of $f$. For every $f\in C([0,1],\R)$, and every $r>0$,
we also set:
$$\omega_f(r)=\sup_{s,t\in[0,1],|t-s|\leq r} |f(t)-f(s)|.$$

We fix $x\geq 0$ and unless otherwise indicated, 
we argue under $\ov\P^n_x$.
 Let $F$ be a bounded 
Lipschitz function on $C([0,1],\R)^2$. We have to prove that
$$\lim_{n\to\infty} \ov\E^n_x[F(C^{(n)},V^{(n)})]=E[F(\ov e^0,\ov Z^0)].$$

We may and will assume that $0\leq F\leq 1$ and that the
Lipschitz constant of $F$ is less than $1$. As in Section 1, for every $r>0$, we denote by
$(\eg,Z^r)$ a Brownian snake with initial point $r$ and
we let $(\ov\eg^r,\ov Z^r)$ be distributed as $(\eg,Z^r)$
conditioned on the event
$$\Big\{\inf_{0\leq t\leq 1} Z^r(t) >0\}\,,$$
which has positive probability. We know that
\be
\label{convdis}
(\ov\eg^r,\ov Z^r)\build{\la}_{r\to 0}^{\rm (d)} (\ov\eg^0,\ov Z^0).
\ee

\begin{lemma}
\label{key-invar}
Let $0<c'<c''$. Then,
$$\lim_{p\to\infty}\;\sup_{c'p^{1/4}\leq y\leq c''p^{1/4}}
|\ov\E^p_{y}[F(C^{(p)},V^{(p)})]-E[F(\ov \eg^{\kappa y/p^{1/4}},\ov Z^{\kappa
y/p^{1/4}})]| =0$$
where $\kappa={1\over \rho}({\sigma\over 2})^{1/2}$.
\end{lemma}

\proof First note that the law of the infimum of a linear Brownian snake
driven by a normalized Brownian excursion $\eg$ has no atoms: See Lemma 2.1 in
\cite{LGW} for the case of an unnormalized Brownian excursion
$e=(e(t),t\geq 0)$ under the It\^o measure, and then use the fact that,
for every $\varepsilon >0$, the law of $(\eg(t),0\leq t\leq 1-\varepsilon)$
is absolutely continuous with respect to that of $(e(t),0\leq t\leq 1-\varepsilon)$. It follows that
the law of $(\ov\eg^r,\ov Z^r)$ depends continuously on $r$. It then suffices to
show that if $(y_p)$ is a sequence such that $c'p^{1/4}\leq y_p\leq c''p^{1/4}$
and $p^{-1/4}y_p\la r$, then,
\be
\label{key-invar1}
\ov\E^p_{y_p}[F(C^{(p)},V^{(p)})]\build{\la}_{p\to\infty}^{} E[F(\ov\eg^{\kappa
r},\ov Z^{\kappa r})].
\ee
However, Theorem \ref{invar0} implies that
$$\E^p_{y_p}[F(C^{(p)},V^{(p)})\,{\bf 1}_{\{\un U
>0\}}]\build{\la}_{p\to\infty}^{} E[F(\eg,Z^{\kappa r})
{\bf 1}_{\{\un Z^{\kappa r}>0\}}]
$$
where $\un Z^{\kappa r}=\inf\{Z^{\kappa r}(t)\,,\,0\leq t\leq 1\}$ (we use
the fact that $P(\un Z^{\kappa r}=0)=0$, as noted above). The desired result
(\ref{key-invar1}) readily follows. \cq

\smallskip
Let $b>0$. We will prove that for $n$
sufficiently large,
$$|\ov\E^n_x[F(C^{(n)},V^{(n)})]-E[F(\ov\eg^{0},\ov Z^{0})]|\leq 12 b,$$
which is enough to get the desired convergence. 

By (\ref{convdis}), we can choose
$\varepsilon\in(0,b\wedge {1\over 100})$ small enough so that
\be
\label{invar5}
|E[F(\ov\eg^r,\ov Z^r)]-E[F(\ov\eg^0,\ov Z^0)]|<b
\ee
for every $0<r\leq 2\varepsilon$. By taking $\varepsilon$ smaller
if necessary, we can also assume that, for every $r\in(0,1]$,
\begin{eqnarray}
\label{invar6}
&&E[(3\varepsilon\sup_{0\leq t\leq 1} \ov\eg^r(t))\wedge 1]\leq b\,,\nonumber\\
\noalign{\smallskip}
&&E[\omega_{\ov\eg^r}(6\varepsilon)\wedge 1]\leq b\,,\nonumber\\
\noalign{\smallskip}
&&E[(3\varepsilon\sup_{0\leq t\leq 1} \ov Z^r(t))\wedge 1]\leq b\,,\nonumber\\
\noalign{\smallskip}
&&E[\omega_{\ov Z^r}(6\varepsilon)\wedge 1]\leq b\;.
\end{eqnarray}
For $\delta,\alpha>0$, denote by
$\Gamma_n=
\Gamma^{\varepsilon,\alpha,\delta}_n$ the event
$$\Gamma_n=\Big\{\inf_{t\in[\delta,1-\delta]}V^{(n)}(t)\geq \alpha\;,\;
\sup_{t\in[0,3\delta]\cup[1-3\delta,1]}(C^{(n)}(t)+V^{(n)}(t))\leq
\varepsilon\Big\}.$$
By Proposition \ref{keytech},
we can fix $\delta,\alpha\in(0,\varepsilon)$ such that, for every
$n$ sufficiently large, 
$$\ov\P^n_x(\Gamma_n)>1-b.$$
On the event $\Gamma_n$, we have, for every $t\in[2n\delta,2n(1-\delta)]$,
\be
\label{invar99}
V(t)\geq \rho\,\Big({2\over \sigma}\Big)^{1/2}\,\alpha\,n^{1/4}.
\ee

Set $\ov\alpha=\rho(2/\sigma)^{1/2}\alpha$. The next step of the proof is
to apply Corollary \ref{spatial}
with $a=\ov\alpha\,n^{1/4}$, assuming that $n$ is large enough so that
$a>x$. Let us introduce the relevant notation. We denote by
$v^n_1,\ldots,v^n_{M_n}$ the exit vertices from $(-\infty,\ov\alpha\,n^{1/4})^2$,
listed in lexicographical order. As in Section 5, we can then consider the
spatial trees $(\t^{[v^n_1]},\ov U^{[v^n_1]}),\ldots,(\t^{[v^n_{M_n}]},\ov
U^{[v^n_{M_n}]})$. The contour functions of these spatial trees may
be obtained in the following way. Set
\ba
&&k^n_1=\inf\{k\in\N: V(k)\geq \ov\alpha\,n^{1/4}\}\\
&&\ell^n_1=\inf\{k\geq k^n_1:C(k+1)<C(k^n_1)\}
\ea
and, by induction on $i$,
\ba
&&k^n_{i+1}=\inf\{k>\ell^n_i: V(k)\geq \ov\alpha\,n^{1/4}\}\\
&&\ell^n_{i+1}=\inf\{k\geq k^n_{i+1}:C(k+1)<C(k^n_{i+1})\}.
\ea
Then $k^n_i\leq \ell^n_i<\infty$ iff $i\leq M_n$. The contour function 
of $\t^{[v^n_i]}$ is 
$$(C(k^n_i+t)-C(k^n_i),0\leq t\leq  \ell^n_i- k^n_i)$$
and the spatial contour function of $(\t^{[v^n_i]},\ov U^{[v^n_i]})$
is
$$(V(k^n_i+t),0\leq t\leq  \ell^n_i- k^n_i).$$
Note in particular that $U_{v^n_i}=V(k^n_i)=V(\ell^n_i)$, and that 
$\ell^n_i-k^n_i=2(|\t^{[v^n_i]}|-1)$.

By construction, for every integer $k\in[0,2n]\backslash
\cup_{i=1}^{M_n}[k^n_i,\ell^n_i]$, we have $V(k)<\ov\alpha\,n^{1/4}$. Also note
that $\ell^n_i+1<k^n_{i+1}$ for every $i\in\{1,\ldots,M_n-1\}$.

Using (\ref{invar99}), we then see that on the event
$\Gamma_n$ all integer points of
$[2n\delta,2n(1-\delta)]$ must be contained in a single interval
$[k^n_i,\ell^n_i]$. Hence, if
$$E_n:=\{\exists i\in\{1,\ldots,M_n\}: \ell^n_i-k^n_i>2(1-3\delta)n\}$$
we have $\Gamma_n\subset E_n$ if $n$ is sufficiently large, and so
$$\ov\P^n_x(E_n)\geq \ov\P^n_x(\Gamma_n)>1-b.$$

On the event $E_n$, we denote by $i_n$ the unique integer $i\in\{1,\ldots,M_n\}$
such that $\ell^n_i-k^n_i>2(1-3\delta)n$. We also set
\ba
&&m_n=|\t^{[v^n_{i_n}]}|-1\,,\\
&&Y_n=U_{v^n_{i_n}}\,.
\ea
Then Corollary \ref{spatial} implies that under the measure $\ov\P^n_x(\cdot\mid
E_n)$, conditionally on the $\sigma$-field
$$\g_n:=\sigma\Big((\t^{\ov\alpha n^{1/4}},U^{\ov\alpha n^{1/4}}),M_n,
(|\t^{[v^n_i]}|,1\leq i\leq M_n)\Big)\,,$$
the spatial tree $(\t^{[v^n_{i_n}]},\ov U^{[v^n_{i_n}]})$ has distribution
$\ov \P^{m_n}_{Y_n}$. Note that $E_n\in\g_n$ and that $Y_n$ and $m_n$
are $\g_n$-measurable.

\begin{lemma}
\label{exitpoint}
The law of $n^{-1/4}Y_n$ under the measure $\ov\P^n_x(\cdot\mid
E_n)$ converges as $n\to\infty$ to the Dirac measure at $\ov\alpha$.
\end{lemma}

\proof For every $v\in\u^*$, write $\check v$ for the father of $v$. 
By construction, we have on $E_n$,
\ba
&&Y_n=U_{v^n_{i_n}}\geq \ov\alpha\, n^{1/4}\\
&&U_{\check v^n_{i_n}}<\ov\alpha\,n^{1/4}.
\ea
To get the statement of the lemma, it thus suffices to verify that, for every
$r>0$,
\be
\label{contin}
\ov\P^n_x\Big(\sup_{v\in\t\backslash\{\varnothing\}}
{|U_v-U_{\check v}|\over n^{1/4}}>r\Big) \build{\la}_{n\to\infty}^{} 0.
\ee
If $\ov \P^n_x$ is replaced by $\P^n_x$, or by $\Q^n$, (\ref{contin})
becomes a straightforward consequence of (\ref{4mom}) (it can also be read 
from Theorem \ref{invar0}). We can then use Lemma
\ref{reroot-equiv} once again to see that (\ref{contin}) also holds when $\ov \P^n_x$
is replaced by $\ov\Q^n$. Finally, the same arguments as in
the end of the proof of Proposition \ref{positive} give 
(\ref{contin}) in the form stated above.\cq

\medskip
On the event $E_n$, we define for $0\leq t\leq 1$,
\begin{eqnarray}
\label{invar1}
&&\wt C^{(n)}(t)={\sigma\over 2}\,
{C(k^n_{i_n}+2m_nt)-C(k^n_{i_n})\over m_n^{1/2}}\nonumber\\
&&\wt V^{(n)}(t)={1\over \rho}\,\Big({\sigma\over 2}\Big)^{1/2}\,
{V(k^n_{i_n}+2m_nt)\over m_n^{1/4}}\;.
\end{eqnarray}
Note that $\wt C^{(n)}$ and $\wt V^{(n)}$ are rescaled versions of the contour
function and the spatial contour function of the spatial tree
$(\t^{[v^n_{i_n}]},\ov U^{[v^n_{i_n}]})$. On the event $E_n^c$ we take
$\wt C^{(n)}(t)=\wt V^{(n)}(t) =0$ for every $0\leq t\leq 1$.

By the remarks preceding Lemma \ref{exitpoint}, we have, for any nonnegative
measurable function $F$ on $C([0,1],\R)^2$,
\be
\label{invar2}
\ov\E^n_x[1_{E_n}F(\wt C^{(n)},\wt V^{(n)})]
=\ov\E^n_x[1_{E_n}\,\ov\E^n_x[F(\wt C^{(n)},\wt V^{(n)})\mid\g_n]]
=\ov\E^n_x[1_{E_n}\,\ov\E^p_{Y_n}[F(C^{(p)},V^{(p)})]_{p=m_n}].
\ee

\medskip
We will be able to combine Lemma \ref{key-invar} and Lemma \ref{exitpoint} in
order to study the right-hand side of (\ref{invar2}). Still
we need to explain why $(C^{(n)},V^{(n)})$ is close to 
$(\wt C^{(n)},\wt V^{(n)})$ under $\ov\P^n_x$, in a suitable sense. 
Recall the definition of the
event
$\Gamma_n\subset E_n$, and the fact that
$\ov \P^n_x(\Gamma_n^c)<b$. Simple estimates show that
for all $n$ sufficiently large  
we have on $\Gamma_n$, for every $0\leq t\leq 1$,
$$|\wt C^{(n)}(t)-C^{(n)}(t)|
\leq \varepsilon + (1-{m_n^{1/2}\over n^{1/2}})\sup_{0\leq s\leq 1} \wt
C^{(n)}(s) +\omega_{\wt C^{(n)}}(6 \delta).$$
In the previous inequality, we used the bounds
$$n\geq m_n\geq (1-3\delta)n\quad,\quad k^n_{i_n}<3\delta n$$
that hold on $\Gamma_n$ for $n$ large. Since 
$$1-{m_n^{1/2}\over n^{1/2}}\leq 1-{m_n\over n} \leq 3\delta\,,$$
we finally get on $\Gamma_n$
\be 
\label{invar3}
\sup_{0\leq t\leq 1} |\wt C^{(n)}(t)-C^{(n)}(t)|
\leq \varepsilon+3\delta \sup_{0\leq t\leq 1} \wt
C^{(n)}(t)+\omega_{\wt C^{(n)}}(6 \delta).
\ee

Similarly, we have on $\Gamma_n$
\be
\label{invar4}
\sup_{0\leq t\leq 1} |\wt V^{(n)}(t)-V^{(n)}(t)|
\leq \varepsilon+3\delta \sup_{0\leq t\leq 1} \wt
V^{(n)}(t)+\omega_{\wt V^{(n)}}(6 \delta).
\ee

Let us now complete the proof. By Lemma \ref{key-invar}, if $p$ is 
sufficiently large, we have
$$\sup_{{\ov\alpha\over 2}p^{1/4}\leq y\leq 2\ov\alpha p^{1/4}}
|\ov\E^p_{y}[F(C^{(p)},V^{(p)})]-E[F(\ov\eg^{\kappa y/p^{1/4}},\ov Z^{\kappa
y/p^{1/4}})]|<b\;.$$
Since $2\kappa \ov\alpha=2\alpha<2\varepsilon$, we can combine this with 
(\ref{invar5}) to get
$$\sup_{{\ov\alpha\over 2}p^{1/4}\leq y\leq 2\ov\alpha p^{1/4}}
|\ov\E^p_{y}[F(C^{(p)},V^{(p)})]-E[F(\ov\eg^{0},\ov Z^{0})]|<2b\;.$$
Now recall (\ref{invar2}), the fact that $m_n\geq (1-3\delta)n$ on $E_n$
and Lemma \ref{exitpoint}. It follows that for $n$ sufficiently large,
$$|\ov\E^n_x[1_{E_n}F(\wt C^{(n)},\wt V^{(n)})]-\ov\P^n_x(E_n)\,
E[F(\ov\eg^{0},\ov Z^{0})]|<3b\;.$$
Since $\ov\P^n_x(E_n)>1-b$, this implies
\be 
\label{invar7}
|\ov\E^n_x[F(\wt C^{(n)},\wt V^{(n)})]-
E[F(\ov\eg^{0},\ov Z^{0})]|<4b\;.
\ee
Furthermore, from the bounds (\ref{invar3}) and (\ref{invar4}), we have
\ba
\ov\E^n_x[1_{\Gamma_n}|F(\wt C^{(n)},\wt V^{(n)})-F(C^{(n)},V^{(n)})|]
\!\!&\leq&\!\! 2\varepsilon+ \ov\E^n_x\Big[\Big(3\delta \sup_{0\leq t\leq 1}
\wt C^{(n)}(t)\Big)\wedge 1+\omega_{\wt C^{(n)}}(6\delta)\wedge 1\Big]\\
&&+\ov\E^n_x\Big[\Big(3\delta \sup_{0\leq t\leq 1}
\wt V^{(n)}(t)\Big)\wedge 1+\omega_{\wt V^{(n)}}(6\delta)\wedge 1\Big]
\ea
At this point, we can again use (\ref{invar2}), Lemma \ref{key-invar} 
and Lemma \ref{exitpoint}
to see that the right-hand side is bounded above for $n$ sufficiently large 
by
\ba
&&b+2\varepsilon+\sup_{0< r\leq 1} E\Big[\Big(3\varepsilon\sup_{0\leq t\leq 1}
\ov\eg^r(t)\Big)\wedge 1 + \omega_{\ov\eg^r}(6\varepsilon)\wedge 1\Big]\\
&&\quad +\sup_{0< r\leq 1} E\Big[\Big(3\varepsilon\sup_{0\leq t\leq 1}
\ov Z^r(t)\Big)\wedge 1 + \omega_{\ov Z^r}(6\varepsilon)\wedge 1\Big]
\ea
From (\ref{invar6}), the latter quantity is bounded above by $7b$.
Since $\ov\P^n_x(\Gamma_n)>1-b$, this gives the bound
$$\ov\E^n_x[|F(\wt C^{(n)},\wt V^{(n)})-F(C^{(n)},V^{(n)})|]\leq 8b.$$
Combining this bound with (\ref{invar7}) leads to
$$|\ov\E^n_x[F(C^{(n)},V^{(n)})]-E[F(\ov\eg^{0},\ov Z^{0})]|\leq 12 b,$$
which completes the proof of Theorem \ref{invar}. \cq

\section{An application to random quadrangulations}

In this section, we apply Theorem \ref{invar} to give a short derivation
of some asymptotics for random quadrangulations which were obtained 
in \cite{CS}. Let us briefly recall the main definitions, following Section
2 of \cite{CS}. A {\it planar map} is a proper embedding, without edge
crossings, of a connected graph in the plane. Loops and multiple edges
are a priori allowed. A planar map is {\it rooted} if there is a
distinguished edge on the border of the infinite face, which is called the
root edge. By convention, the root edge is oriented counterclockwise, and its
origin is called the root vertex. The set of vertices will always be equipped
with the graph distance: If $a$ and $a'$ are two vertices, $d(a,a')$ is the
minimal number of edges on a path from $a$ to $a'$. Two rooted planar maps
are said to be equivalent if there exists a homemorphism of the plane that sends
one map onto the other one and preserves the root edges. 

A planar map is a {\it quadrangulation} if all faces have degree $4$. A
quadrangulation contains no loop but may contain multiple edges. 
For every integer $n\geq 2$, we denote
by $\q_n$ the set of all (equivalent classes of) quadrangulations with
$n$ faces. Then $\q_n$ is a finite set, whose cardinality was computed 
by Tutte \cite{Tu}:
$$|\q_n|={2\over n+2}\,{3^n\over n+1}\,\left(\!\!
\begin{array}{r}
2n\\
n
\end{array}
\right)
.
$$

The relations between planar maps and the present work come from
a basic result (Cori-Vauquelin \cite{CV}, Schaeffer \cite{Sch})
connecting quadrangulations with the so-called well-labelled trees.
Let us call labelled tree any
spatial tree $(\t,U)$ such that $U_\varnothing=1$, $U_v\in\Z$
for every $v\in\t$ and $|U_v-U_{\check v}|\leq 1$ for every 
$v\in\t\backslash\{\varnothing\}$ (recall that $\check v$ is the father of $v$).
The tree is said to be well-labelled if in addition $U_v\geq 1$
for every $v\in\t$. We denote by $\T_n$ the collection
of all labelled trees with $n+1$ vertices, and by $\T^0_n$ the
collection of all well-labelled trees with $n+1$ vertices.

\begin{theorem} 
\label{bijection}
There exists a bijection $\Phi_n$ from $\q_n$ onto $\T^0_n$,
which enjoys the following additional property. Let $q\in\q_n$
and $(\t,U)=\Phi_n(q)$. Then, if ${\cal V}_q$ denotes the set
of vertices of $q$, and $a_0\in{\cal V}_q$ is the root vertex of $q$, we
have for every integer $k\geq 1$:
$$|\{a\in\v_q:d(a_0,a)=k\}|=|\{v\in\t:U_v=k\}|.$$
\end{theorem} 

See Section 3 of \cite{CS} for a detailed proof. If the quadrangulation $q$
and the well-labelled tree $(\t,U)$ are related by the bijection of the theorem,
there is a one-to-one correspondence between vertices $a$ of 
the quadrangulation $q$ other than the
root vertex and vertices $v$ of the tree $\t$, in such a way that the 
distance $d(a_0,a)$ from the root coincides with the label (or spatial position)
$U_v$. This explains the final formula of the theorem.

Before stating the main asymptotic result, let us introduce the 
relevant notation. If $q$ is a rooted quadrangulation, the radius $r(q)$ 
is the maximal distance between the root vertex $a_0$ and another vertex $a$. The
profile
$\lambda_q$ is the integer-valued measure on $\N$ defined by
$$\lambda_q(k)=|\{a\in\v_q:d(a_0,a)=k\}|.$$
Note that $r(q)$ is just the supremum of the support of $\lambda_q$. 
It is also convenient to introduce the rescaled profile. If $q\in\q_n$,
this is the probability measure on $\R_+$ defined by
$$\lambda^{(n)}_q(A)={1\over n+1}\,\lambda_q(n^{1/4}A)$$
for any Borel subset $A$ of $\R_+$.

\begin{theorem}
\label{quadrangulations}
\begin{description}
\item{\bf (i)} The law of $n^{-1/4}r(q)$ under the uniform probability measure
on $\q_n$ converges as $n\to\infty$ to the law 
of the variable
$$\Big({8\over 9}\Big)^{1/4}\;(\sup_{0\leq s\leq 1}  Z^0(s) -\inf_{0\leq s\leq 1}  Z^0(s)).$$
\item{\bf (ii)} The law of the random measure $\lambda^{(n)}_q$
under the uniform probability measure
on $\q_n$ converges as $n\to\infty$ to the law 
of the random probability measure ${\cal I}$ defined by
$$\langle {\cal I},g\rangle =\int_0^1 dr\;g\Big(\Big({8\over 9}\Big)^{1/4} (
Z^0(r) -\inf_{0\leq s\leq 1}  Z^0(s))\Big).$$
\item{\bf (iii)} The law of the rescaled distance $n^{-1/4}d(a_0,a)$ from
a vertex $a$ chosen uniformly at random among all vertices of $q$ to the
root vertex $a_0$, under the uniform probability measure
on $\q_n$, converges as $n\to\infty$ to the law 
of the random variable
$$\Big({8\over 9}\Big)^{1/4}\;(\sup_{0\leq s\leq 1} Z^0(s)).$$
\end{description}
\end{theorem}

\noindent{\bf Remarks.} (a) Part (i) of the theorem is in Corollary 3 of 
\cite{CS} (which also gives the convergence of moments). Part (ii) is
Corollary 4 of \cite{CS}. Part (iii) is not stated in \cite{CS},
but as we will see it is a straightforward consequence of (ii).

(b) We could also have given the various limits in Theorem 
\ref{quadrangulations} in terms of the random measure known as
(one-dimensional) ISE. Up to the trivial multiplicative constant
$(8/9)^{1/4}$, the limit in (i) is the length of the support of ISE,
the limit in (iii) is the supremum of this support, and the random measure
$\cal I$ appearing in (ii) is ISE itself shifted by the minimum
of its support. As is justified precisely in \cite{LGW}, this shifting
is equivalent to conditioning ISE to be suported on the positive half-line. 

(c) Detailed information about the limiting laws in (i) and (iii)
can be found in Delmas \cite{Delmas} and in the recent preprint
Bousquet-M\'elou \cite{Bo}.

\smallskip
\noindent{\bf Proof.} We apply the results of the preceding sections
taking
$\mu(k)=2^{-k-1}$ for $k\in\Z_+$
and letting $\gamma$ be the uniform probability measure on $\{-1,0,1\}$:
$\gamma(-1)=\gamma(0)=\gamma(1)={1\over 3}$.
Note that we have then $\sigma^2=2,\,\rho^2=2/3$ and thus $\kappa={1\over\rho}
\,({\sigma\over 2})^{1/2}=({9\over 8})^{1/4}$.

With the preceding choice of $\mu$ and $\gamma$, one immediately verifies
that $\P^n_1$ is the uniform probability measure on $\T_n$,
and $\ov\P^n_1$ is the uniform probability measure on $\T^0_n$. The
various assertions of Theorem \ref{quadrangulations} can then be obtained
by combining Theorem \ref{bijection} with Theorem \ref{invar}.

To begin with, Theorem \ref{bijection} entails that the law of $r(q)$
under the uniform probability measure
on $\q_n$ coincides with the law of $\sup\{U_v:v\in\t\}$
under $\ov\P^n_1$. Since by construction, if $(\t,U)\in \T_n$,
$$\sup\{U_v:v\in\t\}=\sup\{V(t):t\in[0,2n]\}$$
Theorem \ref{invar} readily implies that the law of $n^{-1/4}\sup\{U_v:v\in\t\}$
under $\q_n$ converges to the law of 
$$\Big({8\over 9}\Big)^{1/4}\;(\sup_{0\leq s\leq 1}  \ov Z^0(s)).$$
From the ``Verwaat transformation'' connecting the conditioned Brownian
snake and the unconditioned one (cf Section 1),
this is the same as the limit in (i).

Let us turn to (ii). By Theorem \ref{bijection}, the law of 
$\lambda^{(n)}_q$ under the uniform probability measure
on $\q_n$ coincides with the law under $\ov\P^n_1$ of the
random measure ${\cal I}_n$ defined by 
$$\langle{\cal I}_n,g\rangle={1\over n+1}\sum_{v\in\t} g(n^{-1/4}U_v).$$
In view of our asymptotics, we may replace ${\cal I}_n$ by ${\cal I}'_n$
defined by
$$\langle{\cal I}'_n,g\rangle={1\over n}\sum_{v\in\t\backslash\{\varnothing\}}
g(n^{-1/4}U_v).$$
Now, from the definition of the contour function $C$ and of the spatial contour
function $V$, it is elementary to verify that
we have also
$$\langle{\cal I}'_n,g\rangle={1\over 2n}\int_0^{2n} dt\,g\Big({V([t]_C)
\over n^{1/4}}\Big)$$
where if $t\in[k,k+1)$ we set $[t]_C=k$ if $C(k)\geq C(t)$ and 
$[t]_C=k+1$ otherwise. In this form, and using the fact that 
$|[t]_C-t|\leq 1$, we deduce from Theorem \ref{invar} that the
law of ${\cal I}'_n$ under $\ov\P^n_1$ converges to the law 
of the random measure ${\cal I}'$ defined by
$$\langle{\cal I}',g\rangle =\int_0^1 dr\,g\Big(\Big({8\over 9}\Big)^{1/4} \ov Z^0(r)\Big).$$
Again the Verwaat transformation shows that this is the same as the limit in
(ii).

Finally, let $X_n$ be distributed as $n^{-1/4}d(a_0,a)$ when the
quadrangulation $q$
is uniform over ${\cal Q}_n$ and $a$ is uniform over the set of vertices of 
$q$ other than the root vertex $a_0$, and let $g$ be bounded and
continuous on $\R_+$. Then, 
$$E[g(X_n)]={1\over |\q_n|} \, \sum_{q\in\q_n} \int \lambda^{(n)}_q(dx)\,g(x) .$$
From (ii), this converges towards
$$E\left(\int_0^1 dr\;g\Big(\Big({8\over 9}\Big)^{1/4} (Z^0(r) -\inf_{0\leq s\leq 1}
Z^0(s))\Big)\right).$$ 
Now, by the invariance property of the Brownian snake
under uniform re-rooting (see e.g. Theorem 
2.3 in \cite{LGW}), the latter quantity is equal to
$$E\left(g\Big(-\Big({8\over 9}\Big)^{1/4}\inf_{0\leq s\leq 1}
Z^0(s)\Big)\right)=E\left(g\Big(\Big({8\over 9}\Big)^{1/4}\sup_{0\leq s\leq 1}
Z^0(s)\Big)\right),$$ 
by symmetry. This completes the proof. \cq

\smallskip
Let us conclude with some remarks. 
The cardinality of $\T_n$ is $3^n$ times the cardinality of the set
of rooted ordered trees with $n+1$ vertices, which is the Catalan number of
order $n$:
$$|\T_n|={3^n\over n+1}\left(\!\!
\begin{array}{r}
2n\\
n
\end{array}
\right).$$
Comparing with the formula for $|\q_n|=|\T^0_n|$, we see that
$$\P^n_1(\un U>0)={|\T^0_n|\over |\T_n|}={2\over n+2}$$
(cf Theorem 2 in \cite{CS} for a combinatorial explanation). This is of course
consistent with the estimates of Proposition \ref{positive}.

The proofs in \cite{CS} are based 
on a form of Theorem \ref{invar0} (which allows one to deal with
labelled trees) and some
delicate combinatorial arguments that are needed to relate
well-labelled trees with labelled trees (the latter are called embedded
trees in \cite{CS}). The originality of our approach 
is thus to apply asymptotics for well-labelled trees, viewed here
as conditioned trees, rather than to use a combinatorial method to
get rid of the conditioning. We expect that this method will have
applications to other types of planar maps, which are also known 
to be in one-to-one correspondence with various classes of discrete trees
(see in particular \cite{BDG3}).

\end{document}